# Tail behaviour of multiple random integrals and $U$-statistics[*]

## Péter Major

*Alfréd Rényi Mathematical Institute of the Hungarian Academy of Sciences*
*e-mail:* `major@renyi.hu`
*url: http://www.renyi.hu/∼major*

**Abstract:** This paper contains sharp estimates about the distribution of multiple random integrals of functions of several variables with respect to a normalized empirical measure, about the distribution of $U$-statistics and multiple Wiener–Itô integrals with respect to a white noise. It also contains good estimates about the supremum of appropriate classes of such integrals or $U$-statistics. The proof of most results is omitted, I have concentrated on the explanation of their content and the picture behind them. I also tried to explain the reason for the investigation of such questions. My goal was to yield such a presentation of the results which a non-expert also can understand, and not only on a formal level.

Received January 2005.

## 1. Formulation of the main problems

To formulate the main problems discussed in this paper first I introduce some notations. Let us have a sequence of independent and identically distributed random variables $\xi_1, \ldots, \xi_n$ on a measurable space $(X, \mathcal{X})$ with distribution $\mu$, and introduce their empirical distribution

$$\mu_n(A) = \frac{1}{n} \#\{j \colon \xi_j \in A, \ 1 \le j \le n\}, \quad A \in \mathcal{X}. \tag{1}$$

Given a measurable function $f(x_1, \ldots, x_k)$ on the product space $(X^k, \mathcal{X}^k)$ let us consider the integral of this function with respect to the $k$-fold direct product of the normalized version $\sqrt{n}(\mu_n - \mu)$ of the empirical measure $\mu_n$, i.e. take the integral

$$J_{n,k}(f) = \frac{n^{k/2}}{k!} \int' f(x_1, \ldots, x_k)(\mu_n(\,dx_1) - \mu(\,dx_1)) \ldots (\mu_n(\,dx_k) - \mu(\,dx_k)),$$

where the prime in $\int'$ means that the diagonals $x_j = x_l$,

$1 \le j < l \le k$, are omitted from the domain of integration. $\tag{2}$

---

[*]This is an original survey paper

I am interested in the following two problems:

*Problem a).* Give a good estimate of the probabilities $P(J_{n,k}(f) > u)$ under some appropriate (and not too restrictive) conditions on the function $f$.

It seems to be natural to omit the diagonals $x_j = x_l$, $j \neq l$, from the domain of integration in the definition of the random integrals $J_{n,k}(f)$ in (2). In the applications I met the estimation of such a version of the integrals was needed.

I shall also discuss the following more general problem:

*Problem b).* Let $f \in \mathcal{F}$ be a nice class of functions on the space $(X^k, \mathcal{X}^k)$. Give a good estimate of the probabilities $P\left(\sup_{f \in \mathcal{F}} J_{n,k}(f) > u\right)$ where $J_{n,k}(f)$ denotes again the random integral of a function $f$ defined in (2).

I met the problems formulated above when I tried to adapt the method of investigation about the limit behaviour of maximum likelihood estimates to more difficult problems, to so-called non-parametric maximum likelihood estimates. An important step in the investigation of maximum likelihood estimates consists of a good approximation of the maximum-likelihood function whose root we are looking for. The Taylor expansion of this function yields a good approximation if its higher order terms are dropped. In an adaptation of this method to more complicated situations the solution of the above mentioned problems a) and b) appear in a natural way. They play a role similar to the estimation of the coefficients of the Taylor expansion in the study of maximum likelihood estimates. Here I do not discuss the details of this approach to non-parametric maximum-likelihood problems. The interested reader may find some further information about it in papers [23] and [24], where such a question is investigated in detail in a special case.

In the above mentioned papers the so-called Kaplan–Meyer method is investigated for the estimation of a distribution function by means of censored data. The solution of problem a) is needed to bound the error of the Kaplan–Meyer estimate for a single argument of the distribution function, and the solution of problem b) helps to bound the difference of this estimate and the real distribution function in the supremum norm. Let me remark that the approach in papers [23] and [24] seems to be applicable under much more general circumstances, but this requires the solution of some hard problems.

I do not know of other authors who dealt directly with the study of random integrals similar to that defined in (2). On the other hand, several authors investigated the behaviour of $U$-statistics, and discussed the next two problems that I describe under the name problem a') and problem b'). To formulate them first I recall the notion of $U$-statistics.

If a sequence of independent and identically distributed random variables $\xi_1, \dots, \xi_n$ is given on a measurable space $(X, \mathcal{X})$ together with a function of $k$ variables $f(x_1, \ldots, x_k)$ on the space $(X^k, \mathcal{X}^k)$, $n \geq k$, then the expression

$$I_{n,k}(f) = \frac{1}{k!} \sum_{\substack{1 \leq j_s \leq n, \ s=1,\ldots,k \\ j_s \neq j_{s'} \text{ if } s \neq s'}} f\left(\xi_{j_1}, \ldots, \xi_{j_k}\right). \tag{3}$$

is called a $U$-statistic of order $k$ with kernel function $f$. Now I formulate the following two problems.

*Problem a').* Give a good estimate of the probabilities $P(n^{-k/2} I_{n,k}(f) > u)$ under some appropriate (and not too restrictive) conditions on the function $f$.

*Problem b').* Let $\mathcal{F}$ be a nice class of functions on the space $(X^k, \mathcal{X}^k)$. Give a good estimate of the probabilities $P\left(\sup_{f \in \mathcal{F}} n^{-k/2} I_{n,k}(f) > u\right)$ where $I_{n,k}(f)$ denotes again the $U$-statistic with kernel function $f$ defined in (3).

Problems a) and b) are closely related to problems a') and b'), but the detailed description of their relation demands some hard work. The main difference between these two pairs of problems is that integration with respect to a power of the measure $\mu_n - \mu$ in formula (2) means some kind of normalization, while the definition of the $U$-statistics in (3) contains no normalization. Moreover, there is no simple way to introduce some good normalization in $U$-statistics. This has the consequence that in problems a) and b) a good estimate can be given for a much larger class of functions than in problems a') and b'). Hence the original pair of problems seems to be more useful in several possible applications.

Both the integrals $J_{n,k}(f)$ defined in (2) and the $U$-statistics $I_{n,k}(f)$ defined in (3) are non-linear functionals of independent random variables, and the main difficulty arises in their study because of this non-linearity. On the other hand, the normalized empirical measure $\sqrt{n}(\mu_n - \mu)$ is close to a Gaussian field for a large sample size $n$. Moreover, as we shall see, $U$-statistics with a large sample size behave similarly to multiple Gaussian integrals. This suggests that the study of multiple Gaussian integrals may help a lot in the solution of our problems. To investigate them first I recall the definition of white noise that we shall need later.

**Definition of a white noise with some reference measure.** *Let us have a $\sigma$-finite measure $\mu$ on a measurable space $(X, \mathcal{X})$. A white noise with reference measure $\mu$ is a Gaussian field $\mu_W = \{\mu_W(A) : A \in \mathcal{X}, \ \mu(A) < \infty\}$, i.e. a set of jointly Gaussian random variables indexed by the above sets $A$, which satisfies the relations $E\mu_W(A) = 0$ and $E\mu_W(A)\mu_W(B) = \mu(A \cap B)$.*

*Remark:* In the definition of a white noise one also mentions the property $\mu_W(A \cup B) = \mu_W(A) + \mu_W(B)$ with probability 1 if $A \cap B = \emptyset$, and $\mu(A) < \infty$, $\mu(B) < \infty$. This could be omitted from the definition, because it follows from the remaining properties of white noises. Indeed, simple calculation shows that $E(\mu_W(A \cup B) - \mu_W(A) - \mu_W(B))^2 = 0$ if $A \cap B = \emptyset$, hence $\mu_W(A \cup B) - \mu_W(A) - \mu_W(B) = 0$ with probability 1 in this case. It also can be observed that if some sets $A_1, \ldots, A_k \in \mathcal{X}$, $\mu(A_j) < \infty$, $1 \leq j \leq k$, are disjoint, then the random variables $\mu_W(A_j)$, $1 \leq j \leq k$, are independent because of the uncorrelatedness of these jointly Gaussian random variables.

It is not difficult to see that for an arbitrary reference measure $\mu$ on a space $(X, \mathcal{X})$ a white noise $\mu_W$ with this reference measure really exists. This follows simply from Kolmogorov's fundamental theorem, by which if the finite dimensional distributions of a random field are prescribed in a consistent way, then there exists a random field with these finite dimensional distributions.

If a white noise $\mu_W$ with a $\sigma$-finite reference measure $\mu$ is given on some measurable space $(X, \mathcal{X})$ together with a function $f(x_1, \ldots, x_k)$ on $(X^k, \mathcal{X}^k)$ such that

$$\sigma^2 = \int f^2(x_1, \ldots, x_k) \mu(dx_1) \ldots \mu(dx_k) < \infty, \tag{4}$$

then the multiple Wiener–Itô integral of the function $f$ with respect to a white noise $\mu_W$ with reference measure $\mu$ can be defined, (see e.g. [14] or [17]). It will be denoted by

$$Z_{\mu,k}(f) = \int f(x_1, \ldots, x_k) \mu_W(dx_1) \ldots \mu_W(dx_k). \tag{5}$$

Here we shall not need a detailed discussion of Wiener–Itô integrals, it will be enough to recall the idea of their definition.

Let us have a measurable space $(X, \mathcal{X})$ together with a non-atomic $\sigma$-finite measure $\mu$ on it. (Wiener–Itô integrals are defined only with respect to a white noise $\mu_W$ with a non-atomic reference measure $\mu$.) We call a function $f$ on $(X^k, \mathcal{X}^k)$ elementary if there exists a finite partition $A_1, \ldots, A_M$, $1 \leq M < \infty$, of the set $X$ (i.e. $A_j \cap A_{j'} = \emptyset$ if $j \neq j'$ and $\bigcup_{j=1}^{M} A_j = X$) such that $\mu(A_j) < \infty$ for all $1 \leq j \leq M$, and the function $f$ satisfies the properties

$$f(x_1, \ldots, x_k) = c(j_1, \ldots, j_k) \quad \text{if } x_1 \in A_{j_1}, \ldots, x_k \in A_{j_k},$$
$$1 \leq j_s \leq M, \ \ 1 \leq s \leq k,$$
$$\text{and } c(j_1, \ldots, j_k) = 0 \quad \text{if } j_s = j_{s'} \ \text{ for some } 1 \leq s < s' \leq k \tag{6}$$

with some real numbers $c(j_1, \ldots, j_k)$, $1 \leq j_s \leq M$, $1 \leq s \leq k$, i.e. the function $f$ is constant on all $k$-dimensional rectangles $A_{j_1} \times \cdots \times A_{j_k}$, and it equals zero on those rectangles which have two sides which agree. (More precisely, we allow the exception $\mu(A_M) = \infty$, but in the case $\mu(A_M) = \infty$ we demand in formula (6) that $c(j_1, \ldots, j_k) = 0$ if one of the arguments $j_s$, $1 \leq s \leq k$ equals $M$. In this case we omit from the sum in the next formula (7) those indices $(j_1, \ldots, j_k)$ for which one of the coordinates of this vector equals $M$.)

The Wiener-Itô integral of the elementary function $f(x_1, \ldots, x_k)$ in formula (6) with respect to a white noise $\mu_W$ with the (non-atomic) reference measure

$\mu$ is defined by the formula

$$
\begin{aligned}
Z_{\mu,k}(f) &= \int f(x_1,\ldots,x_k)\mu_W(dx_1)\ldots\mu_W(dx_k) \\
&= \sum_{1\le j_s\le M,\ 1\le s\le k} c(j_1,\ldots,j_k)\mu_W(A_{j_1})\cdots\mu_W(A_{j_k}).
\end{aligned} \tag{7}
$$

Then the definition of Wiener–Itô integral can be extended to a general function satisfying relation (4) by means of an $L_2$-isomorphism. The details of this extension will be not discussed here.

Let me remark that the condition $c(j_1,\ldots,j_k)=0$ if $j_s=j_{s'}$ for some $1\le s < s'\le k$ in the definition of an elementary functions can be interpreted so that, similarly to the definition of the random integral $J_{n,k}(f)$ in (2), the diagonals are omitted from the domain of integration of a Wiener–Itô integral $Z_\mu(f)$.

The investigation of Wiener–Itô integrals is simpler than that of random integrals $J_{n,k}(f)$ defined in (2) or of $U$-statistics introduced in (3) because of the Gaussian property of the underlying white noise. Beside this, the study of Wiener–Itô integrals may help in understanding what kind of estimates can be expected in the solution of problems a) and b) or a') and b') and also in finding the proofs. Hence it is useful to consider the following two problems.

*Problem a'').* Give a good estimate of the probabilities $P(Z_{\mu,k}(f)>u)$ under some appropriate (and not too restrictive) conditions on the function $f$ and measure $\mu$.

*Problem b'').* Let $\mathcal{F}$ be a nice class of functions on the space $(X^k,\mathcal{X}^k)$. Give a good estimate of the probabilities $P\left(\sup_{f\in\mathcal{F}} Z_{\mu,k}(f)>u\right)$ where $Z_{\mu,k}(f)$ denotes again a Wiener–Itô integral with function $f$ and white noise with reference measure $\mu$.

In this paper the above problems will be discussed. Such estimates will be presented which depend on some basic characteristics of the random expressions $J_{n,k}(f)$, $I_{n,k}(f)$ or $Z_{\mu,k}(f)$. They will depend mainly on the $L_2$ and $L_\infty$-norm of the function $f$ taking part in the definition of the above quantities. (The $L_2$-norm of the function $f$ is closely related to the variance of the random variables we consider.) The proof of the estimates is related to some other problems interesting in themselves. My main goal was to explain the results and ideas behind them. I put emphasis on the explanation of the picture that can help understanding them, and the details of almost all proofs are omitted. A detailed explanation together with the proofs can be found in my Lecture Note [22].

This paper consists of 9 sections. The first four sections contain the results about problems a), a') and a'') together with some other statements which may explain better their background. Section 5 contains the main ideas of their proof. In Section 6 problems b), b') and b'') are discussed together with some related questions. The main ideas of the proofs of the results in Section 6 which contain many unpleasant technical details are discussed in Sections 7 and 8. In Section 9 Talagrand's theory about concentration inequalities is considered together with some new results and open questions.

## 2. The discussion of some large deviation results

First we restrict our attention to problems a), a′) and a″), i.e. to the case when the distribution of the random integral or $U$-statistic of one function is estimated. These problems are much simpler in the special case $k = 1$. But they are not trivial even in this case. A discussion of some large deviation results may help to understand them better. I recall some large deviation results, but not in their most general form. Actually these results will not be needed later, they are interesting for the sake of some orientation.

**Theorem 2.1 (Large deviation theorem about partial sums of independent and identically distributed random variables).** *Let $\xi_1, \xi_2, \ldots$, be a sequence of independent and identically distributed random variables such that $E\xi_1 = 0$, $Ee^{t\xi_1} < \infty$ with some $t > 0$. Let us define the partial sums $S_n = \sum\limits_{j=1}^{n} \xi_j$, $n = 1, 2, \ldots$. Then the relation*

$$\lim_{n \to \infty} \frac{1}{n} \log P(S_n \geq nu) = -\rho(u) \qquad \text{for all } u > 0 \tag{8}$$

*holds with the function $\rho(u)$ defined by the formula $\rho(u) = \sup\limits_{t} \left( tu - \log Ee^{t\xi_1} \right)$.*

*The function $\rho(\cdot)$ in formula (8) has the following properties: $\rho(u) > 0$ for all $u > 0$, and it is a monotone increasing function, there is some number $0 < A \leq \infty$ with a number $A$ depending on the distribution of the function $\xi_1$ such that $\rho(u) < \infty$ for $0 \leq u \leq A$, and the asymptotic relation $\rho(u) = \frac{\sigma^2 u^2}{2} + O(u^3)$ holds for small $u > 0$, where $\sigma^2 = E\xi_1^2$ is the variance of $\xi_1$.*

The above theorem states that for all $\varepsilon > 0$ the inequality $P(S_n > nu) \leq e^{-n(\rho(u)-\varepsilon)}$ holds if $n \geq n(u, \varepsilon)$, and this estimate is essentially sharp. Actually, in nice cases, when the equation $\rho(u) = \sup\limits_{t} \left( tu - \log Ee^{t\xi_1} \right)$ has a solution in $t$, the above inequality also holds with $\varepsilon = 0$ for all $n \geq 1$. The function $\rho(u)$ in the exponent of the above large deviation estimate strongly depends on the distribution of $\xi_1$. It is the so-called Legendre transform of $\log Ee^{t\xi_1}$, of the logarithm of the moment generating function of $\xi_1$, and its values in an arbitrary interval determine the distribution of $\xi_1$. On the other hand, the estimate (8) for small $u > 0$ shows some resemblance to the bound suggested by the central limit theorem. Indeed, for small $u > 0$ it yields the upper bound $e^{-n\sigma^2 u^2/2 + nO(u^3)}$, while the central limit theorem would suggest the estimate $e^{-n\sigma^2 u^2/2}$. (Let us recall that the standard normal distribution function $\Phi(u)$ satisfies the inequality $\left( \frac{1}{u} - \frac{1}{u^3} \right) \frac{e^{-u^2/2}}{\sqrt{2\pi}} < 1 - \Phi(u) < \frac{1}{u} \frac{e^{-u^2/2}}{\sqrt{2\pi}}$ for all $u > 0$, hence for large $u$ it is natural to bound it by $e^{-u^2/2}$.)

The next result I mention, Bernstein's inequality, (see e.g. [5], 1.3.2 Bernstein's inequality) has a closer relation to the problems discussed in this paper. It gives a good upper bound on the distribution of sums of independent, bounded random variables with expectation zero. It is important that this estimate is universal, the constants it contains do not depend on the properties of the random variables we consider.

**Theorem 2.2 (Bernstein's inequality).** *Let $X_1, \ldots, X_n$ be independent random variables, $P(|X_j| \leq 1) = 1$, $EX_j = 0$, $1 \leq j \leq n$. Put $\sigma_j^2 = EX_j^2$, $1 \leq j \leq n$, $S_n = \sum\limits_{j=1}^{n} X_j$ and $V_n^2 = Var\, S_n = \sum\limits_{j=1}^{n} \sigma_j^2$. Then*

$$P(S_n > u) \leq \exp\left\{ -\frac{u^2}{2V_n^2\left(1 + \frac{u}{3V_n^2}\right)} \right\} \quad \text{for all } u > 0. \tag{9}$$

Let us take a closer look on the content of Theorem 2.2. Estimate (9) yields a bound of different form if the first term is dominating in the sum $1 + \frac{u}{3V_n^2}$ in the denominator of the fraction in this expression and if the second term is dominating in it. If we fix some constant $C > 0$, then formula (9) yields that $P(S_n > u) \leq e^{-Bu^2/2V_n^2}$ with some constant $B = B(C)$ for $0 \leq u \leq CV_n^2$. If, moreover $0 \leq u \leq \varepsilon V_n^2$ with some small $\varepsilon > 0$, then the estimate $P(S_n > u) \leq e^{-(1-K\varepsilon)u^2/2V_n^2}$ holds with a universal constant $K > 0$. This means that in the case $0 < u \leq CV_n^2$ the tail behaviour of the distribution of $F(u) = P(S_n > u)$ can be bounded by the distribution $G(u) = P(\text{const.}\, V_n \eta > u)$ where $\eta$ is a standard normal random variable, and $V_n^2$ is the variance of the partial sum $S_n$. If $0 \leq u \leq \varepsilon V_n^2$ with a small $\varepsilon > 0$, then it also can be bounded by $P((1 - K\varepsilon)V_n \eta > u)$ with some universal constant $K > 0$.

In the case $u \gg V_n^2$ formula (9) yields a different type of estimate. In this case we get that $P(S_n > u) < e^{-(3-\varepsilon)u/2}$ with a small $\varepsilon > 0$, and this seems to be a rather weak estimate. In particular, it does not depend on the variance $V_n^2$ of $S_n$. In the degenerate case $V_n = 0$ when $P(S_n > u) = 0$, estimate (9) yields a strictly positive upper bound for $P(S_n > u)$. One would like to get such an improvement of Bernstein's inequality which gives a better bound in the case $u \gg V_n^2$. Bennett's inequality (see e.g. [28], Appendix B, 4 Bennett's inequality) satisfies this requirement.

**Theorem 2.3 (Bennett's inequality).** *Let $X_1, \ldots, X_n$ be independent random variables, $P(|X_j| \leq 1) = 1$, $EX_j = 0$, $1 \leq j \leq n$. Put $\sigma_j^2 = EX_j^2$, $1 \leq j \leq n$, $S_n = \sum\limits_{j=1}^{n} X_j$ and $V_n^2 = Var\, S_n = \sum\limits_{j=1}^{n} \sigma_j^2$. Then*

$$P(S_n > u) \leq \exp\left\{ -V_n^2\left[ \left(1 + \frac{u}{V_n^2}\right)\log\left(1 + \frac{u}{V_n^2}\right) - \frac{u}{V_n^2} \right] \right\} \quad \text{for all } u > 0. \tag{10}$$

*As a consequence, for all $\varepsilon > 0$ there exists some $B = B(\varepsilon) > 0$ such that*

$$P(S_n > u) \leq \exp\left\{-(1-\varepsilon)u\log\frac{u}{V_n^2}\right\} \quad \text{if } u > BV_n^2, \tag{11}$$

*and there exists some positive constant $K > 0$ such that*

$$P(S_n > u) \leq \exp\left\{-Ku\log\frac{u}{V_n^2}\right\} \quad \text{if } u > 2V_n^2. \tag{12}$$

Estimates (11) or (12) yield a slight improvement of Bernstein's inequality in the case $u \geq KV_n^2$ with a sufficiently large $K > 0$. On the other hand, even this estimate is much weaker than the estimate suggested by a formal application of the central limit theorem. The question arises whether they are sharp or can be improved. The next example shows that inequalities (11) or (12) in Bennett's inequality are essentially sharp. If no additional restrictions are imposed, then at most the universal constants can be improved in them. Even a sum of independent, bounded and identically distributed random variables can be constructed which satisfies a lower bound similar to the upper bounds in formulas (11) and (12), only with possibly different constants.

**Example 2.4.** *Let us fix some positive integer $n$, real numbers $u$ and $\sigma^2$ such that $0 < \sigma^2 \leq \frac{1}{8}$, $n > 3u \geq 6$ and $u > 4n\sigma^2$. Put $V_n^2 = n\sigma^2$ and take a sequence of independent, identically distributed random variables $X_1, \ldots, X_n$ such that $P(X_j = 1) = P(X_j = -1) = \frac{\sigma^2}{2}$, and $P(X_j = 0) = 1 - \sigma^2$. Put $S_n = \sum_{j=1}^{n} X_j$. Then $ES_n = 0$, $Var S_n = V_n^2$, and*

$$P(S_n \geq u) > \exp\left\{-Bu\log\frac{u}{V_n^2}\right\}$$

*with some appropriate (universal) constant $B > 0$.*

*Remark:* The estimate of Example 2.4 or of relations (11) and (12) is well comparable with the tail distribution of a Poisson distributed random variable with parameter $\lambda = \text{const.}\, n\sigma^2 \geq 1$ at level $u \geq 2\lambda$. Some calculation shows that a Poisson distributed random variable $\zeta_\lambda$ with parameter $\lambda > 1$ satisfies the inequality $e^{-C_1 u\log(u/\lambda)} \leq P(\zeta_\lambda - E\zeta_\lambda > u) \leq P(\zeta_\lambda > u) \leq P(\zeta_\lambda - E\zeta_\lambda > \frac{u}{2}) \leq e^{-C_2 u\log(u/\lambda)}$ with some appropriate constants $0 < C_1 < C_2 < \infty$ for all $u > 2\lambda$, and $E\zeta_\lambda = Var\,\zeta_\lambda = \lambda$. This estimate is similar to the above mentioned relations.

Example 2.4 is proved in Example 3.2 of my Lecture Note [22], but here I present a simpler proof.

*Proof of the statement of Example 2.4.* Let us fix an integer $u$ such that $n > 3u$ and $u > 4n\sigma^2$. Let $B = B(u)$ denote the event that among the random variables $X_j$, $1 \leq j \leq n$, there are exactly $3u$ terms with values $+1$ or $-1$, and all other random variables $X_j$ equal zero. Let us also define the event $A = A(u) \subset B(u)$ which holds if $2u$ random variables $X_j$ are equal to 1, $u$ random variables $X_j$ are equal to $-1$, and all remaining random variables $X_j$, $1 \leq j \leq n$, are equal to zero. Clearly, $P(S_n \geq u) \geq P(A) = P(B)P(A|B)$. On the other hand, $P(B) = \binom{n}{3u} (\sigma^2)^{3u} (1 - \sigma^2)^{n-3u} \geq \left(\frac{n}{3u}\right)^{3u} (\sigma^2)^{3u} e^{-4n\sigma^2} = e^{-3u \log(3u/n\sigma^2) - 4n\sigma^2}$. Here we exploited that because of the condition $\sigma^2 \leq \frac{1}{8}$ we have $1 - \sigma^2 \geq e^{-4\sigma^2}$. Beside this, $u \leq 4n\sigma^2$, and $P(B) \geq e^{-3u \log(3u/n\sigma^2) - u} \geq e^{-B_1 u \log(u/n\sigma^2)}$ with some appropriate $B_1 > 0$ under our assumptions.

Let us consider a set of $3u$ elements, and choose a random subset of it by taking all elements of this set with probability $1/2$ to this random subset independently of each other. I claim that the conditional probability $P(A|B)$ equals the probability that this random subset has $2u$ elements. Indeed, even the conditional probability of the event $A$ under the condition that for a prescribed set of indices $J \subset \{1, \ldots, n\}$ with exactly $3u$ elements we have $X_j = \pm 1$ if $j \in J$ and $X_j = 0$ if $j \notin J$ equals the probability of the event that the above defined random subset has $2u$ elements. This is so, because under this condition the random variables $X_j$ take the value $+1$ with probability $1/2$ for all $j \in J$ independently of each other. Hence $P(A|B) = \binom{3u}{2u} 2^{-3u} \geq e^{-Cu} \geq e^{-B_2 u \log(u/n\sigma^2)}$ with some appropriate constants $C > 0$ and $B_2 > 0$ under our conditions, since $\frac{u}{n\sigma^2} \geq 4$ in this case. The estimates given for $P(B)$ and $P(A|B)$ imply the statement of Example 2.4.

Bernstein's inequality provides a solution to problems a) and a$'$) in the case $k = 1$ under some conditions. Because of the normalization (multiplication by $n^{-1/2}$ in these problems) it yields an estimate with the choice $\bar{u} = \sqrt{n}u$. Observe that $J_{n,1}(f) = \frac{1}{\sqrt{n}} \sum_{j=1}^{n} (f(\xi_j) - Ef(\xi_j))$ for $k = 1$ in the definition (2). In problem a) it gives a good bound on $P(J_{n,1}(f) > u)$ for a function $f$ such that $|f(x)| \leq \frac{1}{2}$ for all $x \in X$ with the choice $X_j = f(\xi_j) - Ef(\xi_j)$, $1 \leq j \leq n$, and $\bar{u} = \sqrt{n}u$. In problem a$'$) it gives a good bound on $P(n^{-1/2}I_{n,1}(f) > u)$ under the condition $|f(x)| \leq 1$ for all $x \in X$, and $Ef(\xi_1) = 0$ with the choice $X_j = f(\xi_j)$, $1 \leq j \leq n$, and $\bar{u} = \sqrt{n}u$. This means that in the case $0 \leq u \leq C\sqrt{n}\sigma^2$ the bounds $P(J_{n,1}(f) > u) \leq e^{-Ku^2/2\sigma^2}$ and $P(n^{-1/2}I_{n,1}(f) > u) \leq e^{-Ku^2/2\sigma^2}$ hold with $\sigma^2 = \text{Var}\, f(\xi_1)$ and some constant $K = K(C)$ depending on the number $C$ if the above conditions are imposed in problem a) or a$'$). If $0 \leq u \leq \varepsilon \sqrt{n}\sigma^2$ with some small $\varepsilon > 0$, then the above constant $K$ can be chosen very close to the number 1.

The above results can be interpreted so that that in the case $0 \leq u \leq \text{const.} \sqrt{n}\sigma^2$ and a bounded function $f$ an estimate suggested by the central limit theorem holds for problem a), only an additional constant multiplier may appear in the exponent. A similar statement holds in problem a$'$), only here the additional condition $Ef(\xi_j) = 0$ has to be imposed. On the other hand, the situation is quite different if $u \gg \sqrt{n}\sigma^2$. In this case Bernstein's inequality yields only a very weak estimate. Bennett's inequality gives a slight improvement. It yields the inequality $P(J_{n,1}(f) > u) \leq e^{-Bu\sqrt{n}\log(u/\sqrt{n}\sigma^2)}$ with an appropriate constant $B > 0$ if $|f(x)| \leq \frac{1}{2}$ for all $x \in X$, $u \geq 2\sqrt{n}\sigma^2$, and $\sigma^2 = \operatorname{Var} f(\xi_1)$. The estimate $P(n^{-1/2}I_{n,1}(f) > u) \leq e^{-Bu\sqrt{n}\log(u/\sqrt{n}\sigma^2)}$ holds with an appropriate $B > 0$ if $|f(x)| \leq 1$ for all $x \in X$, $Ef(\xi_1) = 0$, $\operatorname{Var} f(\xi_1) = \sigma^2$, and $u \geq 2\sqrt{n}\sigma^2$. These estimates are much weaker than the bound suggested by a formal application of the central limit theorem. On the other hand, as Example 2.4 shows, no better estimate can be expected in this case. Moreover, the proof of this example gives some insight why a different type of estimate appears in the cases $u \leq \sqrt{n}\sigma^2$ and $u \gg \sqrt{n}\sigma^2$ for problems a) and a').

In the proof of Example 2.4 a 'bad' irregular event $A$ was defined such that if it holds, then the sum of the random variables considered in this example is sufficiently large. Generally, the probability of such an event is very small, but if the variance of the random variables is very small, (in problems a) and a') this is the case if $\sigma^2 \ll un^{-1/2}$) then such 'bad' irregular events can be defined whose probabilities are not negligible.

Problems a) and a') will also be considered for $k \geq 2$, and this will be called the multivariate case. The results we get for the solution of problems a) and a') in the multivariate case is very similar to the results described above. To understand them first some problems have to be discussed. In particular, the answer for the following two questions has to be understood:

*Question a).* In the solution of problem a') in the case $k = 1$ the condition $Ef(\xi_1) = 0$ was imposed, and this means some kind of normalization. What condition corresponds to it in the multivariate case? This question leads to the definition of degenerate $U$-statistics and to the so-called Hoeffding's decomposition of $U$-statistics to a sum of degenerate $U$-statistics.

*Question b).* The discussion of problems a) and a') was based on the central limit theorem. What kind of limit theorems can take its place in the multivariate case? What kind of limit theorems hold for $U$-statistics $I_{n,k}(f)$ or multiple random integrals $J_{n,k}(f)$ defined in (2)? The limit appearing in these problems can be expressed by means of multiple Wiener–Itô integrals in a natural way.

In the next section the two above questions will be discussed.

## 3. On some problems about $U$-statistics and random integrals

### 3.1. *The normalization of $U$-statistics*

In the case $k = 1$ problem a') means the estimation of sums of independent and identically distributed random variables. In this case a good estimate was obtained under the condition $Ef(\xi_1) = 0$.

In the multivariate case $k \geq 2$ a stronger normalization property has to be imposed to get good estimates about the distribution of $U$-statistics. In this case it has to be assumed that the conditional expectations of the terms $f(\xi_{j_1}, \ldots, \xi_{j_k})$ of the $U$-statistic under the condition that the value of all but one arguments takes a prescribed value equals zero. This property is formulated in a more explicit way in the following definition of degenerate $U$-statistics.

**Definition of degenerate $U$-statistics.** *Let us consider the $U$-statistic $I_{n,k}(f)$ of order $k$ defined in formula (3) with kernel function $f(x_1, \ldots, x_k)$ and a sequence of independent and identically distributed random variables $\xi_1, \ldots, \xi_n$. It is a degenerate $U$-statistic if its kernel function satisfies the relation*

$$E f(\xi_1, \ldots, \xi_k | \xi_1 = x_1, \ldots, \xi_{j-1} = x_{j-1}, \xi_{j+1} = x_{j+1}, \ldots, \xi_k = x_k) = 0$$
$$\text{for all } 1 \leq j \leq k \text{ and } x_s \in X, \ s \in \{1, \ldots, k\} \setminus \{j\}. \tag{13}$$

The definition of degenerate $U$-statistics is closely related to the notion of canonical functions described below.

**Definition of canonical functions.** *A function $f(x_1, \ldots, x_k)$ taking values in the $k$-fold product $(X^k, \mathcal{X}^k)$ of a measurable space $(X, \mathcal{X})$ is called canonical with respect to a probability measure $\mu$ on $(X, \mathcal{X})$ if*

$$\int f(x_1, \ldots, x_{j-1}, u, x_{j+1}, \ldots, x_k) \mu(\,du) = 0$$
$$\text{for all } \ 1 \leq j \leq k \ \text{ and } \ x_s \in X, \ s \in \{1, \ldots, k\} \setminus \{j\}. \tag{14}$$

It is clear that a $U$-statistic $I_{n,k}(f)$ is degenerate if and only if its kernel function $f$ is canonical with respect to the distribution $\mu$ of the random variables $\xi_1, \ldots, \xi_n$ appearing in the definition of the $U$-statistic.

Given a function $f$ and a probability measure $\mu$, this function can be written as a sum of canonical functions (with different sets of arguments) with respect to the measure $\mu$, and this enables us to decompose a $U$-statistic as a linear combination of degenerate $U$-statistics. This is the content of Hoeffding's decomposition of $U$-statistics described below. To formulate it first I introduce some notations.

Consider the $k$-fold product $(X^k, \mathcal{X}^k, \mu^k)$ of a measure space $(X, \mathcal{X}, \mu)$ with some probability measure $\mu$, and define for all integrable functions $f(x_1, \ldots, x_k)$ and indices $1 \leq j \leq k$ the projection $P_j f$ of the function $f$ to its $j$-th coordinate as

$$P_j f(x_1, \ldots, x_{j-1}, x_{j+1}, \ldots, x_k) = \int f(x_1, \ldots, x_k) \mu(\,dx_j), \quad 1 \leq j \leq k. \tag{15}$$

In some investigations it may be useful to rewrite formula (15) by means of conditional expectations in an equivalent form as

$$P_j f(x_1, \ldots, x_{j-1}, x_{j+1}, \ldots, x_k)$$
$$= E(f(\xi_1, \ldots, \xi_k) | \xi_1 = x_1, \ldots, \xi_{j-1} = x_{j-1}, \xi_{j+1} = x_{j+1}, \ldots, \xi_k = x_k),$$

where $\xi_1, \ldots, \xi_k$ are independent random variables with distribution $\mu$.

Let us also define the operators $Q_j = I - P_j$ as $Q_j f = f - P_j f$ on the space of integrable functions on $(X^k, \mathcal{X}^k, \mu^k)$, $1 \leq j \leq k$. In the definition (15) $P_j f$ is a function not depending on the coordinate $x_j$, but in the definition of $Q_j$ we introduce the fictive coordinate $x_j$ to make the expression $Q_j f = f - P_j f$ meaningful. The following result holds.

**Theorem 3.1 (Hoeffding's decomposition of $U$-statistics).** *Let an integrable function $f(x_1, \ldots, x_k)$ be given on the $k$-fold product space $(X^k, \mathcal{X}^k, \mu^k)$ of a space $(X, \mathcal{X}, \mu)$ with a probability measure $\mu$. It has the decomposition*

$$f = \sum_{V \subset \{1, \ldots, k\}} f_V, \quad with$$

$$f_V(x_j, \, j \in V) = \left( \prod_{j \in \{1, \ldots, k\} \setminus V} P_j \prod_{j \in V} Q_j \right) f(x_1, \ldots, x_k) \quad (16)$$

*such that all functions $f_V$, $V \subset \{1, \ldots, k\}$, in (16) are canonical with respect to the probability measure $\mu$, and they depend on the $|V|$ arguments $x_j$, $j \in V$.*

*Let $\xi_1, \ldots, \xi_n$ be a sequence of independent, $\mu$ distributed random variables, and consider the $U$-statistics $I_{n,k}(f)$ and $I_{n,|V|}(f_V)$ corresponding to the kernel functions $f$, $f_V$ defined in (16) and random variables $\xi_1, \ldots, \xi_n$. Then*

$$I_{n,k}(f) = \sum_{V \subset \{1, \ldots, k\}} (n - |V|)(n - |V| - 1) \cdots (n - k + 1) \frac{|V|!}{k!} I_{n,|V|}(f_V) \quad (17)$$

*is a representation of $I_{n,k}(f)$ as a sum of degenerate $U$-statistics, where $|V|$ denotes the cardinality of the set $V$. (The product $(n - |V|)(n - |V| - 1) \cdots (n - k + 1)$ is defined as 1 if $V = \{1, \ldots, k\}$, i.e. $|V| = k$.) This representation is called the Hoeffding decomposition of $I_{n,k}(f)$.*

Hoeffding's decomposition was originally proved in paper [13]. It may be interesting also to mention its generalization in [32].

I omit the proof of Theorem 3.1, although it is fairly simple. I only try to briefly explain that the construction of Hoeffding's decomposition is natural. Let me recall that a random variable can be decomposed as a sum of a random variable with expectation zero plus a constant, and the random variable with expectation zero in this decomposition is defined by taking out from the original random variable its expectation. To introduce such a transformation which turns to zero not only the expectation of the transformed random variable, but also its conditional expectation with respect to some condition it is natural to take out from the original random variable its conditional expectation. Since the operators $P_j$ defined in (15) are closely related to the conditional expectations appearing in the definition of degenerate $U$-statistics, the above consideration makes natural to write the identity $f = \prod_{j=1}^{k}(P_j + Q_j)f = \sum_{V \subset \{1, \ldots, k\}} f_V$ with the functions defined in (16). (In the justification of the last formula some properties of the operators $P_j$ and $Q_j$ have to be exploited.)

It is clear that $EI_{n,k}(f) = 0$ for a degenerate $U$-statistic. Also the inequality

$$E\left(I_{n,k}(f)\right)^2 \leq \frac{n^k}{k!}\sigma^2 \quad \text{with } \sigma^2 = \int f^2(x_1,\ldots,x_k)\mu(\,dx_1)\ldots\mu(\,dx_k) \quad (18)$$

holds if $I_{n,k}(f)$ is a degenerate $U$-statistic. The measure $\mu$ in (18) is the distribution of the random variables $\xi_j$ taking part in the definition of the $U$-statistic. Moreover, $\lim_{n\to\infty} n^{-k}E\left(I_{n,k}(f)\right)^2 = \frac{\sigma^2}{k!}$ if the kernel function $f$ is a symmetric function of its arguments, i.e. $f(x_1,\ldots,x_k) = f(x_{\pi(1)},\ldots,x_{\pi(k)})$ for all permutations $\pi = (\pi(1),\ldots,\pi(k))$ of the set $\{1,\ldots,k\}$.

Relation (18) can be proved by means of the observation that

$$Ef(\xi_{j_1},\ldots,\xi_{j_k})f(\xi_{j'_1},\ldots,\xi_{j'_k}) = 0$$

if $\{j_1,\ldots,j_k\} \neq \{j'_1,\ldots,j'_k\}$, and $f$ is a canonical function with respect to the distribution $\mu$ of the random variables $\xi_j$. On the other hand,

$$|Ef(\xi_{j_1},\ldots,\xi_{j_k})f(\xi_{j'_1},\ldots,\xi_{j'_k})| \leq \int f^2(x_1,\ldots,x_k)\mu(\,dx_1)\ldots\mu(\,dx_k)$$

by the Schwarz inequality if $\{j_1,\ldots,j_k\} = \{j'_1,\ldots,j'_k\}$, i.e. if the sequence of indices $j'_1\ldots,j'_k$ is a permutation of the sequence of indices $j_1,\ldots,j_k$, and there is an identity in this relation if the function $f$ is symmetric. The last formula enables us to check the asymptotic relation given for $E\left(I_{n,k}(f)\right)^2$ after relation (18).

Relation (18) suggests to restrict our attention in the investigation of problem a') to degenerate $U$-statistics, and it also explains why the normalization $n^{-k/2}$ was chosen in it. For degenerate $U$-statistics with this normalization such an upper bound can be expected in problem a') which does not depend on the sample size $n$. The estimation of the distribution of a general $U$-statistic can be reduced to the degenerate case by means of Hoeffding's decomposition (Theorem 3.1).

The random integrals $J_{n,k}(f)$ are defined in (2) by means of integration with respect to the signed measure $\mu_n - \mu$, and this means some sort of normalization. This normalization has the consequence that the distributions of these integrals satisfy a good estimation for rather general kernel functions $f$. Beside this, a random integral $J_{n,k}(f)$ can be written as a sum of $U$-statistics to which the Hoeffding decomposition can be applied. Hence it can be rewritten as a linear combination of degenerate $U$-statistics. In the next result I describe the representation of $J_{n,k}(f)$ we get in such a way. It shows that the implicit normalization caused by integration with respect to $\mu_n - \mu$ has a serious cancellation effect. This enables us to get a good solution for problem a) or b) if we have a good solution for problem a') or b'). Unfortunately, the proof of this result demands rather unpleasant calculations. Hence here I omit the proof. It can be found in [19] or in Theorem 9.4 of [22].

**Theorem 3.2.** *Let us have a non-atomic measure $\mu$ on a measurable space $(X, \mathcal{X})$ together with a sequence of independent, $\mu$-distributed random variables $\xi_1, \ldots, \xi_n$, and take a function $f(x_1, \ldots, x_k)$ of $k$ variables on the space $(X^k, \mathcal{X}^k)$ such that*

$$\int f^2(x_1, \ldots, x_k)\mu(dx_1)\ldots\mu(dx_k) < \infty.$$

*Let us consider the empirical distribution function $\mu_n$ of the sequence $\xi_1, \ldots, \xi_n$ introduced in (1) together with the $k$-fold random integral $J_{n,k}(f)$ of the function $f$ defined in (2). The identity*

$$J_{n,k}(f) = \sum_{V \subset \{1,\ldots,k\}} C(n, k, V) n^{-|V|/2} I_{n,|V|}(f_V), \qquad (19)$$

*holds with the canonical (with respect to the measure $\mu$) functions $f_V(x_j, j \in V)$ defined in (16) and appropriate real numbers $C(n, k, V)$, $V \subset \{1, \ldots, k\}$, where $I_{n,|V|}(f_V)$ is the (degenerate) U-statistic with kernel function $f_V$ and random sequence $\xi_1, \ldots, \xi_n$ defined in (3). The constants $C(n, k, V)$ in (19) satisfy the relations $|C(n, k, V)| \leq C(k)$ with some constant $C(k)$ depending only on the order $k$ of the integral $J_{n,k}(f)$, $\lim\limits_{n \to \infty} C(n, k, V) = C(k, V)$ with some constant $C(k, V) < \infty$ for all $V \subset \{1, \ldots, k\}$, and $C(n, k, \{1, \ldots, k\}) = 1$ for $V = \{1, \ldots, k\}$.*

Let us also remark that the functions $f_V$ defined in (16) satisfy the inequalities

$$\int f_V^2(x_j, j \in V) \prod_{j \in V} \mu(dx_j) \leq \int f^2(x_1, \ldots, x_k)\mu(dx_1)\ldots\mu(dx_k) \qquad (20)$$

and

$$\sup_{x_j, \, j \in V} |f_V(x_j, j \in V)| \leq 2^{|V|} \sup_{x_j, \, 1 \leq j \leq k} |f(x_1, \ldots, x_k)| \qquad (21)$$

for all $V \subset \{1, \ldots, k\}$.

The decomposition of the random integral $J_{n,k}(f)$ in formula (19) is similar to the Hoeffding decomposition of general $U$-statistics presented in Theorem 3.1. The main difference between them is that the coefficients of the normalized degenerate $U$-statistics $n^{-|V|/2} I_{n,|V|}(f_V)$ at the right-hand side of formula (19) can be bounded by a universal constant depending neither on the sample size $n$, nor on the kernel function $f$ of the random integral. This fact has important consequences.

Theorem 3.2 enables us to get good estimates for problem a) if we have such estimates for problem a'). In particular, formulas (18), (19) and (20) yield good bounds on the expectation and variance of the random integral $J_{n,k}(f)$. The inequalities

$$E\left(J_{n,k}(f)\right)^2 \leq C\sigma^2 \quad \text{and} \quad |EJ_{n,k}(f)| \leq C\sigma,$$

$$\text{with} \quad \sigma^2 = \int f^2(x_1, \ldots, x_k)\mu(dx_1)\ldots\mu(dx_k) \qquad (22)$$

hold with some universal constant $C > 0$ depending only on the order of the random integral $J_{n,k}(f)$.

Relation (22) yields such an estimate for the second moment of $J_{n,k}(f)$ as we expect. On the other hand, although it gives a sufficiently good bound on its first moment, it does not state that the expectation of $J_{n,k}(f)$ equals zero. Indeed, formula (19) only gives that $|EJ_{n,k}(f)| = |C(n,k,\emptyset)f_{\emptyset}| \leq C|f_{\emptyset}| = C \left| \int f(x_1, \ldots, x_k) \mu(dx_1) \ldots \mu(dx_k) \right| \leq C\sigma$ with some appropriate constant $C > 0$. The following example shows that $EJ_{n,k}(f)$ need not be always zero. (To understand better why such a situation may appear observe that the random measures $(\mu_n - \mu)(B_1)$ and $(\mu_n - \mu)(B_2)$ are not independent for disjoint sets $B_1$ and $B_2$.)

Let us consider a random integral $J_{n,2}(f)$ of order 2 with an appropriate kernel function $f$. Beside this, choose a sequence of independent random variables $\xi_1, \ldots, \xi_n$ with uniform distribution on the unit interval $[0,1]$ and denote its empirical distribution by $\mu_n$. We shall consider the example where the kernel function $f = f(x, y)$ is the indicator function of the unit square, i.e. $f(x, y) = 1$ if $0 \leq x, y \leq 1$, and $f(x, y) = 0$ otherwise. The random integral $J_{n,2}(f) = n \int_{x \neq y} f(x,y)(\mu_n(dx) - dx)(\mu_n(dy) - dy)$ will be taken, and its expected value $EJ_{n,2}(f)$ will be calculated. By adjusting the diagonal $x = y$ to the domain of integration and taking out the contribution obtained in this way we get that $EJ_{n,2}(f) = nE(\int_0^1 (\mu_n(dx) - \mu(dx))^2 - n^2 \cdot \frac{1}{n^2} = -1$, i.e. the expected value of this random integral is not equal to zero. (The last term is the integral of the function $f(x,y)$ on the diagonal $x = y$ with respect to the product measure $\mu_n \times \mu_n$ which equals $(\mu_n - \mu) \times (\mu_n - \mu)$ on the diagonal.)

Now I turn to the second problem discussed in this section.

### 3.2. Limit theorems for U-statistics and random integrals

The following limit theorem about normalized degenerate U-statistics will be interesting for us.

**Theorem 3.3 (Limit theorem about normalized degenerate U-statistics).** *Let us consider a sequence of degenerate U-statistics $I_{n,k}(f)$ of order $k$, $n = k, k+1, \ldots$, defined in (3) with the help of a kernel function $f(x_1, \ldots, x_k)$ on the $k$-fold product $(X^k, \mathcal{X}^k)$ of a measurable space $(X, \mathcal{X})$, canonical with respect to some non-atomic probability measure $\mu$ on $(X, \mathcal{X})$ and such that $\int f^2(x_1, \ldots, x_k) \mu(dx_1) \ldots \mu(dx_k) < \infty$ together with a sequence of independent and identically distributed random variables $\xi_1, \xi_2, \ldots$ with distribution $\mu$ on $(X, \mathcal{X})$. The sequence of normalized U-statistics $n^{-k/2} I_{n,k}(f)$ converges in distribution, as $n \to \infty$, to the $k$-fold Wiener–Itô integral*

$$\frac{1}{k!} Z_{\mu,k}(f) = \frac{1}{k!} \int f(x_1, \ldots, x_k) \mu_W(dx_1) \ldots \mu_W(dx_k)$$

*with kernel function $f(x_1, \ldots, x_k)$ and a white noise $\mu_W$ with reference measure $\mu$.*

The proof of Theorem 3.3 can be found for instance in [6]. Here I present a heuristic explanation which can be considered as a sketch of proof.

To understand Theorem 3.3 it is useful to rewrite the normalized degenerate $U$-statistics considered in it in the form of multiple random integrals with respect to a normalized empirical measure. The identity

$$n^{-k/2}I_{n,k}(f) = n^{k/2}\int{}' f(x_1,\ldots,x_k)\mu_n(\,dx_1)\ldots\mu_n(\,dx_k) \tag{23}$$

$$= n^{k/2}\int{}' f(x_1,\ldots,x_k)(\mu_n(\,dx_1)-\mu(\,dx_1))\ldots(\mu_n(\,dx_k)-\mu(\,dx_1))$$

holds, where $\mu_n$ is the empirical distribution function of the sequence $\xi_1,\ldots,\xi_n$ defined in (1), and the prime in $\int'$ denotes that the diagonals, i.e. the points $x=(x_1,\ldots,x_k)$ such that $x_j=x_{j'}$ for some pairs of indices $1\le j,j'\le k$, $j\ne j'$ are omitted from the domain of integration. The last identity of formula (23) holds, because in the case of a function $f(x_1,\ldots,x_k)$ canonical with respect to a non-atomic measure $\mu$ we get the same result by integrating with respect to $\mu_n(\,dx_j)$ and with respect to $\mu_n(\,dx_j)-\mu(\,dx_j)$. (The non-atomic property of the measure $\mu$ is needed to guarantee that the integrals with respect to the measure $\mu$ considered in this formula remain zero if the diagonals are omitted from the domain of integration.)

Formula (23) may help to understand Theorem 3.3, because the random fields $n^{1/2}(\mu_n(A)-\mu(A))$, $A\in\mathcal{X}$, converge to a Gaussian field $\nu(A)$, $A\in\mathcal{X}$, as $n\to\infty$, and this suggests a limit similar to the result of Theorem 3.3. But it is not so simple to carry out a limiting procedure leading to the proof of Theorem 3.3 with the help of formula (23). Some problems arise, because the fields $n^{1/2}(\mu_n-\mu)$ converge to a not white noise type Gaussian field. The limit we get is similar to a Wiener bridge on the real line. Hence a relation between Wiener processes and Wiener bridges suggests to write the following version of formula (23). Let $\eta$ be a standard Gaussian random variable, independent of the random sequence $\xi_1,\xi_2,\ldots$. We can write, by exploiting again the canonical property of the function $f$, the identity

$$n^{-k/2}I_{n,k}(f) = n^{k/2}\int{}' \quad f(x_1,\ldots,x_k)(\mu_n(\,dx_1)-\mu(\,dx_1)+\eta\mu(\,dx_1))$$

$$\ldots(\mu_n(\,dx_k)-\mu(\,dx_k)+\eta\mu(\,dx_k)). \tag{24}$$

The random measures $n^{1/2}(\mu_n-\mu+\eta\mu)$ converge to a white noise with reference measure $\mu$, hence a limiting procedure in formula (24) yields Theorem 3.3. Moreover, in the case of elementary functions $f$ the central limit theorem and formula (24) imply the statement of Theorem 3.3 directly. (Elementary functions are defined in formula (6).) After this, Theorem 3.3 can be proved in the general case with the help of the investigation of the $L_2$-contraction property of some operators. I omit the details.

A similar limit theorem holds for random integrals $J_{n,k}(f)$. It can be proved by means of Theorem 3.2 and an adaptation of the above sketched argument for the proof of Theorem 3.3. It states the following result.

**Theorem 3.4. Limit theorem about multiple random integrals $J_{n,k}(f)$.**
*Let us have a sequence of independent and identically distributed random variables $\xi_1, \xi_2, \ldots$ with some non-atomic distribution $\mu$ on a measurable space $(X, \mathcal{X})$ and a function $f(x_1, \ldots, x_k)$ on the k-fold product $(X^k, \mathcal{X}^k)$ of the space $(X, \mathcal{X})$ such that*

$$\int f^2(x_1, \ldots, x_k) \mu(\,dx_1) \ldots \mu(\,dx_k) < \infty.$$

*Let us consider for all $n = 1, 2, \ldots$ the random integrals $J_{n,k}(f)$ of order $k$ defined in formulas (1) and (2) with the help of the empirical distribution $\mu_n$ of the sequence $\xi_1, \ldots, \xi_n$ and the function $f$. The random integrals $J_{n,k}(f)$ converge in distribution, as $n \to \infty$, to the following sum $U(f)$ of multiple Wiener–Itô integrals:*

$$
\begin{aligned}
U(f) &= \sum_{V \subset \{1, \ldots, k\}} \frac{C(k, V)}{|V|!} Z_{\mu, |V|}(f_V) \\
&= \sum_{V \subset \{1, \ldots, k\}} \frac{C(k, V)}{|V|!} \int f_V(x_j, \, j \in V) \prod_{j \in V} \mu_W(dx_j),
\end{aligned}
$$

*where the functions $f_V(x_j \, j \in V)$, $V \subset \{1, \ldots, k\}$, are those functions defined in formula (16) which appear in the Hoeffding decomposition of the function $f(x_1, \ldots, x_k)$, the constants $C(k, V)$ are the limits appearing in the limit relation $\lim\limits_{n \to \infty} C(n, k, V) = C(k, V)$ satisfied by the quantities $C(n, k, V)$ in formula (19), and $\mu_W$ is a white noise with reference measure $\mu$.*

The results of this section suggest that to understand what kind of results can be expected for the solution of problems a) and a′) it is useful to study first their simpler counterpart, problem a″) about multiple Wiener–Itô integrals. They also show that problem a′) is interesting in the case when degenerate $U$-statistics are investigated. The next section contains some results about these problems.

## 4. Estimates on the distribution of random integrals and $U$-statistics

First I formulate the results about the solution of problem a″), about the tail-behaviour of multiple Wiener–Itô integrals.

**Theorem 4.1.** *Let us consider a $\sigma$-finite measure $\mu$ on a measurable space $(X, \mathcal{X})$ together with a white noise $\mu_W$ with reference measure $\mu$. Let us have a real-valued function $f(x_1, \ldots, x_k)$ on the space $(X^k, \mathcal{X}^k)$ which satisfies relation (4) with some $\sigma^2 < \infty$. Take the random integral $Z_{\mu,k}(f)$ introduced in formula (5). It satisfies the inequality*

$$P(|Z_{\mu,k}(f)| > u) \leq C \exp\left\{ -\frac{1}{2} \left( \frac{u}{\sigma} \right)^{2/k} \right\} \quad \text{for all } u > 0 \tag{25}$$

*with an appropriate constant $C = C(k) > 0$ depending only on the multiplicity $k$ of the integral.*

The proof of Theorem 4.1 can be found in my paper [20] together with the following example which shows that it gives a sharp estimate.

**Example 4.2.** *Let us have a $\sigma$-finite measure $\mu$ on some measurable space $(X, \mathcal{X})$ together with a white noise $\mu_W$ on $(X, \mathcal{X})$ with reference measure $\mu$. Let $f_0(x)$ be a real valued function on $(X, \mathcal{X})$ such that $\int f_0(x)^2 \mu(\,dx) = 1$, and take the function $f(x_1, \ldots, x_k) = \sigma f_0(x_1) \cdots f_0(x_k)$ with some number $\sigma > 0$ and the Wiener–Itô integral $Z_{\mu,k}(f)$ introduced in formula (5). Then the relation*

$$\int f(x_1, \ldots, x_k)^2 \,\mu(\,dx_1) \ldots \mu(\,dx_k) = \sigma^2$$

*holds, and the Wiener–Itô integral $Z_{\mu,k}(f)$ satisfies the inequality*

$$P(|Z_{\mu,k}(f)| > u) \geq \frac{\bar{C}}{\left(\frac{u}{\sigma}\right)^{1/k} + 1} \exp\left\{-\frac{1}{2} \left(\frac{u}{\sigma}\right)^{2/k}\right\} \quad \textit{for all } u > 0 \qquad (26)$$

*with some constant $\bar{C} > 0$.*

Let us also remark that a Wiener–Itô integral $Z_{\mu,k}(f)$ defined in (5) with a kernel function $f$ satisfying relation (4) also satisfies the relations $EZ_{\mu,k}(f) = 0$ and $EZ_{\mu,k}(f)^2 \leq k!\sigma^2$ with the number $\sigma^2$ in (4). If the function $f$ is symmetric, i.e. if $f(x_1, \ldots, x_k) = f(x_{\pi(1)}, \ldots, x_{\pi(k)})$ for all permutations $\pi$ of the set $\{1, \ldots, k\}$, then in the last relation identity can be written instead of inequality. Beside this, $Z_{\mu,k}(f) = Z_{\mu,k}(\text{Sym } f)$, where $\text{Sym } f$ denotes the symmetrization of the function $f$, and this means that we can restrict our attention to the Wiener–Itô integrals of symmetric functions without violating the generality. Hence Theorem 4.1 can be interpreted in the following way. The random integral $Z_{\mu,k}(f)$ has expectation zero, its variance is less than or equal to $k!\sigma^2$ under the conditions of this result, and there is identity in this relation if $f$ is a symmetric function. Beside this, the distribution of $Z_{\mu,k}(f)$ satisfies an estimate similar to that of $\sigma\eta^k$, where $\eta$ is a standard normal random variable. The estimate (25) in Theorem 4.1 is not always sharp, but Example 4.2 shows that there are cases when the expression in its exponent cannot be improved.

Let me also remark that the above statement can be formulated in a slightly nicer form if the distribution of $Z_{\mu,k}(f)$ is compared not with that of $\sigma\eta^k$, but with that of $\sigma H_k(\eta)$, where $H_k(x)$ is the $k$-th Hermite polynomial with leading coefficient 1. The identities $EH_k(\eta) = 0$, $EH_k(\eta)^2 = k!$ hold. This means that not only the tail distributions of $Z_{\mu,k}(f)$ and $\sigma H_k(\eta)$ are similar, but in the case of a symmetric function $f$ also their first two moments agree.

In problems a) and a$'$) a slightly weaker but similar estimate holds. In the case of problem a$'$) the following result is valid (see [20]).

**Theorem 4.3.** *Let $\xi_1, \ldots, \xi_n$ be a sequence of independent and identically distributed random variables on a space $(X, \mathcal{X})$ with some distribution $\mu$. Let us consider a function $f(x_1, \ldots, x_k)$ on the space $(X^k, \mathcal{X}^k)$, canonical with respect* *to the measure $\mu$ which satisfies the conditions*

$$\|f\|_\infty \;=\; \sup_{x_j \in X,\, 1 \le j \le k} |f(x_1, \ldots, x_k)| \le 1, \tag{27}$$

$$\|f\|_2^2 \;=\; \int f^2(x_1, \ldots, x_k)\mu(\,dx_1)\ldots\mu(\,dx_k) \le \sigma^2, \tag{28}$$

*with some $0 < \sigma^2 \le 1$ together with the degenerate $U$-statistic $I_{n,k}(f)$ defined in formula (3) with this kernel function $f$. There exist some constants $A = A(k) > 0$ and $B = B(k) > 0$ depending only on the order $k$ of the $U$-statistic $I_{n,k}(f)$ such that*

$$P(k!n^{-k/2}|I_{n,k}(f)| > u) \le A \exp\left\{ -\frac{u^{2/k}}{2\sigma^{2/k}\left(1 + B\left(un^{-k/2}\sigma^{-(k+1)}\right)^{1/k}\right)} \right\} \tag{29}$$

*for all $0 \le u \le n^{k/2}\sigma^{k+1}$.*

*Remark:* Actually, the universal constant $B > 0$ can be chosen independently of the order $k$ of the degenerate $U$-statistic $I_{n,k}(f)$ in inequality (29).

Theorem 4.3 can be considered as a generalization of Bernstein's inequality Theorem 2.2 to the multivariate case in a slightly weaker form when only the sum of independent and identically distributed random variables is considered. Its statement, inequality (29) does not contain an explicit value for the constants $A$ and $B$, which are equal to $A = 2$ and $B = \frac{1}{3}$ in the case of Bernstein's inequality. (The constant $A = 2$ appears, because of the absolute value in the probability at the left-hand side of (29).) There is a formal difference between formula (9) and the statement of formula (29) in the case $k = 1$, because in formula (29) the $U$-statistic $I_{n,k}(f)$ of order $k$ is multiplied by $n^{-k/2}$. Another difference between them is that inequality (29) in Theorem 4.3 is stated under the condition $0 \le u \le n^{k/2}\sigma^{k+1}$, and this restriction has no counterpart in Bernstein's inequality. But, as I shall show, Theorem 4.3 also contains an estimate for $u \ge n^{k/2}\sigma^{k+1}$ in an implicit way, and it can be considered as the multivariate version of Bernstein's inequality.

Bernstein's inequality gives a good estimate only if $0 \le u \le K\sqrt{n}\sigma^2$ with some $K > 0$ (with the normalization of Theorem 4.3, i.e. if the probability

$$P\left(n^{-1/2}\sum_{k=1}^n X_k > u\right)$$

is considered). In the multivariate case a similar picture appears. We get a good estimate for problem a′) suggested by Theorem 4.1 only under the condition $0 \le u \le \text{const.}\, n^{k/2}\sigma^{k+1}$. If $0 < u \le \varepsilon n^{k/2}\sigma^{k+1}$ with a sufficiently small $\varepsilon > 0$, then Theorem 4.3 implies the inequality $P(k!n^{-k/2}|I_{n,k}(f)| > u) \le A \exp\left\{ -\frac{1 - C\varepsilon^{1/k}}{2}\left(\frac{u}{\sigma}\right)^{2/k} \right\}$ with some universal constants $A > 0$ and $C > 0$ depending only on the order $k$ of the $U$-statistic $I_{n,k}(f)$. This means that in this case Theorem 4.3 yields an almost as good estimate as Theorem 4.1 about the distribution of multiple Wiener–Itô integrals. We have seen that Bernstein's inequality has a similar property if the estimate (9) is compared with the central limit theorem in the case $0 < u \le \varepsilon V_n^2$ with a small $\varepsilon > 0$.

To see what kind of estimate Theorem 4.3 yields in the case $u \ge n^{k/2}\sigma^{k+1}$ let us observe that in condition (28) we have an inequality and not an identity. Hence in the case $n^{k/2} \ge u > n^{k/2}\sigma^{k+1}$ relation (29) holds with $\bar\sigma = \left(un^{-k/2}\right)^{1/(k+1)}$, and this yields that

$$
\begin{aligned}
P(k!n^{-k/2}|I_{n,k}(f)| > u) &\le A\exp\left\{-\frac{1}{2(1+B)^{1/k}}\left(\frac{u}{\bar\sigma}\right)^{2/k}\right\} \\
&= Ae^{-(u^2n)^{1/(k+1)}/2(1+B)^{1/k}}.
\end{aligned}
$$

(The inequality $n^{k/2} \ge u$ was imposed to satisfy the condition $0 \le \bar\sigma^2 \le 1$.) If $u > n^{k/2}$, then the probability at the left-hand side of (29) equals zero because of condition (27). It is not difficult to see by means of the above calculation that Theorem 4.3 implies the inequality

$$
P\left(k!n^{-k/2}|I_{n,k}(f)| > u\right) \tag{30}
$$
$$
\le c_1\exp\left\{-\frac{c_2u^{2/k}}{\sigma^{2/k}\left(1 + c_3\left(un^{-k/2}\sigma^{-(k+1)}\right)^{2/k(k+1)}\right)}\right\} \quad \text{for all } u > 0
$$

with some universal constants $c_1$, $c_2$ and $c_3$ depending only on the order $k$ of the $U$-statistic $I_{n,k}(f)$, if the conditions of Theorem 4.3 hold. Inequality (30) holds for all $u \ge 0$. Arcones and Giné formulated and proved this estimate in a slightly different but equivalent form in paper [3] under the name generalized Bernstein's inequality. This result is weaker than Theorem 4.3, since it does not give a good value for the constant $c_2$. The method of paper [3] is based on a symmetrization argument. Symmetrization arguments can be very useful in the study of problems b) and b$'$) formulated in the Introduction, but they cannot supply a proof of Theorem 4.3 with good constants because of some principal reasons.

The following result which can be considered as a solution of problem a) is a fairly simple consequence of Theorem 4.3, Theorem 3.2 and formulas (20) and (21).

**Theorem 4.4.** *Let us take a sequence of independent and identically distributed random variables $\xi_1, \ldots, \xi_n$ on a measurable space $(X, \mathcal{X})$ with a non-atomic distribution $\mu$ on it together with a measurable function $f(x_1, \ldots, x_k)$ on the $k$-fold product $(X^k, \mathcal{X}^k)$ of the space $(X, \mathcal{X})$ with some $k \ge 1$ which satisfies conditions (27) and (28) with some constant $0 < \sigma \le 1$. Then there exist some constants $C = C_k > 0$ and $\alpha = \alpha_k > 0$ such that the random integral $J_{n,k}(f)$ defined by formulas (1) and (2) with this sequence of random variables $\xi_1, \ldots, \xi_n$* *and function f satisfies the inequality*

$$P\left(|J_{n,k}(f)| > u\right) \le C \exp\left\{-\alpha \left(\frac{u}{\sigma}\right)^{2/k}\right\} \quad \text{for all } 0 < u \le n^{k/2}\sigma^{k+1}. \quad (31)$$

Theorem 4.4 provides a slightly weaker estimate on the probability considered in Problem a) than Theorem 4.3 about its counterpart in Problem a'). It does not give an almost optimal constant $\alpha$ in the inequality (31) for $0 \le u \le \varepsilon n^{k/2}\sigma^{k+1}$ with a small $\varepsilon > 0$. On the other hand, this estimate is sharp in that sense that disregarding the value of the universal constant $\alpha$ it cannot be improved. It seems to be appropriate in the solution of the problems about non-parametric maximum likelihood estimates mentioned in the Introduction.

The estimate (31) on the probability $P\left(|J_{n,k}(f)| > u\right)$ can be rewritten, similarly to relation (30), in such a form which holds for all $u > 0$. On the other hand, both Theorem 4.3 and Theorem 4.4 yield a very weak estimate if $u \gg n^{k/2}\sigma^{k+1}$. We met a similar situation in Section 2 when these problems were investigated in the case $k = 1$. It is natural to expect that a generalization of Bennett's inequality holds in the multivariate case $k \ge 2$, and it gives an improvement of estimates (29) and (31) in the case $u \gg n^{k/2}\sigma^{k+1}$ for all $k \ge 1$. I can prove only partial results in this direction which are not sharp in the general case. On the other hand, there is a possibility to give such a generalization of Example 2.4 which shows that the inequalities implied by Theorem 4.3 or 4.4 in the case $u \ge n^{k/2}\sigma^{k+1}$, $k \ge 2$ have only a slight improvement.

The results of Theorems 4.3 and 4.4 imply that in the case $u \le n^{k/2}\sigma^{k+1}$ under the condition of these results the probabilities $P(n^{k/2}|I_{n,k}(f)| > u)$ and $P(|J_{n,k}(f)| > u)$ can be bounded by $P(C\sigma|\eta|^k > u)$ with an appropriate universal constant $C = C(k) > 0$ depending only on the order $k$ of the degenerate U-statistic $I_{n,k}(f)$ or of the multiple random integral $J_{n,k}(f)$, where the random variable $\eta$ has standard normal distribution, and

$$\sigma^2 = \int f^2(x_1, \ldots, x_k)\mu(dx_1)\ldots\mu(dx_k).$$

A generalization of Example 2.4 can be given which shows for all $k \ge 1$ that in the case $u \gg n^{k/2}\sigma^{k+1}$ we can have only a much weaker estimate. I shall present such an example only for $k = 2$, but it can be generalized for all $k \ge 1$. This example is taken from my Lecture Note [22] (Example 8.6). Here I present it without a detailed proof. The proof which exploits the properties of Example 2.4 is not long. But I found more instructive to explain the idea behind this example.

**Example 4.5.** *Let $\xi_1, \ldots, \xi_n$ be a sequence of independent, identically distributed valued random variables taking values in the plane, i.e. in $X = R^2$, such that $\xi_j = (\eta_{j,1}, \eta_{j,2})$, $\eta_{j,1}$ and $\eta_{j,2}$ are independent, $P(\eta_{j,1} = 1) = P(\eta_{j,1} = -1) = \frac{\sigma^2}{8}$, $P(\eta_{j,1} = 0) = 1 - \frac{\sigma^2}{4}$, $P(\eta_{j,2} = 1) = P(\eta_{j,2} = -1) = \frac{1}{2}$ for all $1 \le j \le n$. Let us introduce the function $f(x, y) = f((x_1, x_2), (y_1, y_2)) = x_1 y_2 + x_2 y_1$, $x = (x_1, x_2) \in R^2$, $y = (y_1, y_2) \in R^2$ on $X^2$, and define the U-statistic*

$$I_{n,2}(f) = \sum_{1 \le j, k \le n, \, j \ne k} (\eta_{j,1}\eta_{k,2} + \eta_{k,1}\eta_{j,2}) \quad (32)$$

*of order 2 with the above kernel function $f$ and the sequence of independent random variables $\xi_1, \ldots, \xi_n$. Then $I_{n,2}(f)$ is a degenerate U-statistic. If $u \geq B_1 n \sigma^3$ with some appropriate constant $B_1 > 0$, $B_2^{-1} n \geq u \geq B_2 n^{-2}$ with a sufficiently large fixed number $B_2 > 0$, and $1 \geq \sigma \geq \frac{1}{n}$, then the estimate*

$$P(n^{-1} I_{n,2}(f) > u) \geq \exp\left\{ -B n^{1/3} u^{2/3} \log\left(\frac{u}{n\sigma^3}\right) \right\} \tag{33}$$

*holds with some constant $B > 0$ depending neither on $n$ nor on $\sigma$.*

It is not difficult to see that the U-statistic $I_{n,2}(f)$ introduced in Example 4.5 is a degenerate U-statistic of order two with a kernel function $f$ such that $\sup |f(x,y)| \leq 1$ and $\sigma^2 = \int f^2(x,y)\mu(dx)\mu(dy) = E(2\eta_{j,1}\eta_{j,2})^2 = \sigma^2$. Example 4.5 means that in the case $u \gg n\sigma^3$, (i.e. if $u \gg n^{k/2}\sigma^{k+1}$ with $k = 2$) a much weaker estimate holds than in the case $u \leq n\sigma^3$. Let us fix the numbers $u$ and $n$, and consider the dependence of our estimate on $\sigma$. The estimate $P(n^{-1}|I_{n,2}(f)| > u) \leq e^{-Ku/\sigma} = e^{-Ku^{2/3}n^{1/3}}$ holds if $\sigma = u^{1/3}n^{-1/3}$, and Example 4.5 shows that a rather weak improvement appears if $\sigma \ll u^{1/3}n^{-1/3}$.

To understand why the statement of Example 4.5 holds observe that a small error is made if the condition $j \neq k$ is omitted from the summation in formula (32), and this suggests that the approximation

$$\frac{1}{n} I_{n,2}(f) \sim \frac{2}{n} \left( \sum_{j=1}^n \eta_{j,1} \right) \left( \sum_{j=1}^n \eta_{j,2} \right)$$

causes a negligible error. This fact together with the independence of the sequences $\eta_{j,1}$, $1 \leq j \leq n$, and $\eta_{j,2}$, $1 \leq j \leq n$, imply that

$$
\begin{aligned}
P(n^{-1} I_{n,2}(f) > u) &\sim P\left( \left( \sum_{j=1}^n \eta_{j,1} \right) \left( \sum_{j=1}^n \eta_{j,2} \right) > \frac{nu}{2} \right) \\
&\geq P\left( \sum_{j=1}^n \eta_{j,1} > v_1 \right) P\left( \sum_{j=1}^n \eta_{j,2} > v_2 \right)
\end{aligned}
\tag{34}
$$

with such a choice of numbers $v_1$ and $v_2$ for which $v_1 v_2 = \frac{nu}{2}$.

The first probability at the right-hand side of (34) can be bounded because of the result of Example 2.4 as $P\left( \sum_{j=1}^n \eta_{j,1} > v_1 \right) \geq e^{-Bv_1 \log(4v_1/n\sigma^2)}$ if $v_1 \geq 4n\sigma^2$, and the second probability as $P\left( \sum_{j=1}^n \eta_{j,2} > v_2 \right) \geq C e^{-Kv_2^2/n}$ with some appropriate $C > 0$ and $K > 0$ if $0 \leq v_2 \leq n$. The proof of Example 4.5 can be obtained by means of an appropriate choice of the numbers $v_1$ and $v_2$.

In Theorem 4.1 the distribution of a $k$-fold Wiener–Itô integral $Z_{\mu,k}(f)$ was bounded by the distribution of $\sigma\eta^k$ with a standard normal random variable $\eta$ and an appropriate constant $\sigma$. By Theorems 4.3 and 4.4 a similar, but weaker estimate holds for the distribution of a degenerate $U$-statistic $I_{n,k}(f)$ or random integral $J_{n,k}(f)$. In the next section I briefly explain why such results hold.

There is a method to get a good estimate on the moments of the random variables considered in the above theorems, and they enable us to get a good estimate also on the distribution of the random integrals and $U$-statistics appearing in these theorems. The moments of a $k$-fold Wiener–Itô integral can be bounded by the moments of $\sigma\eta^k$ with an appropriate $\sigma > 0$, and this estimate implies Theorem 4.1. Theorems 4.3 and 4.4 can be proved in a similar way. But we can give a good estimate only on not too high moments of the random variables $I_{n,k}(f)$ and $J_{n,k}(f)$, and this is the reason why we get only a weaker result for their distribution.

*Remark:* My goal was to obtain a good estimate in Problems a) and a$'$) if we have a bound on the $L_2$ and $L_\infty$ norm of the kernel function $f$ in them. A similar problem was considered in Problem a$''$) about Wiener–Itô integrals with the difference that in this case only an $L_2$ bound of the function $f$ is needed. Theorems 4.1, 4.3 and 4.4 provided such a bound, and as Example 4.2 shows these estimates are sharp. On the other hand, if we have some additional information about the kernel function $f$, then more precise estimates can be given which in certain cases yield an essential improvement. Such results were known for $U$-statistics and Wiener–Itô integrals of order $k = 2$, (see [9] and [12]) and quite recently (after the submission of the first version of this work) they were generalized in [1] and [15] to general $k \geq 2$. Moreover, these improvements are useful in the study of some problems. Hence a referee suggested to explain them in the present work. I try to follow his advice by inserting their discussion at the end, in the open problems part of the paper.

## 5. On the proof of the results in Section 4

Theorem 4.1 can be proved by means of the following

**Proposition 5.1.** *Let the conditions of Theorem 4.1 be satisfied for a multiple Wiener–Itô integral $Z_{\mu,k}(f)$ of order $k$. Then, with the notations of Theorem 4.1, the inequality*

$$E\left(|Z_{\mu,k}(f)|\right)^{2M} \leq 1 \cdot 3 \cdot 5 \cdots (2kM - 1)\sigma^{2M} \tag{35}$$

*holds for all $M = 1, 2, \ldots$.*

By the Stirling formula Proposition 5.1 implies that

$$E(|Z_{\mu,k}(f)|)^{2M} \leq \frac{(2kM)!}{2^{kM}(kM)!}\sigma^{2M} \leq A\left(\frac{2}{e}\right)^{kM}(kM)^{kM}\sigma^{2M} \tag{36}$$

for any $A > \sqrt{2}$ if $M \geq M_0 = M_0(A)$, and this estimate is sharp. The following Proposition 5.2 which can be applied in the proof of Theorem 4.3 states a similar, but weaker inequality for the moments of normalized degenerate $U$-statistics.

**Proposition 5.2.** *Let us consider a degenerate $U$-statistic $I_{n,k}(f)$ of order $k$ with sample size $n$ and with a kernel function $f$ satisfying relations (27) and (28) with some $0 < \sigma^2 \le 1$. Fix a positive number $\eta > 0$. There exist some universal constants $A = A(k) > \sqrt{2}$, $C = C(k) > 0$ and $M_0 = M_0(k) \ge 1$ depending only on the order of the $U$-statistic $I_{n,k}(f)$ such that*

$$E \left( n^{-k/2} k! I_{n,k}(f) \right)^{2M} \le A \left( 1 + C\sqrt{\alpha} \right)^{2kM} \left( \frac{2}{e} \right)^{kM} (kM)^{kM} \sigma^{2M} \tag{37}$$

*for all integers $M$ such that $k M_0 \le kM \le \alpha n \sigma^2$.*

*The constant $C = C(k)$ in formula (37) can be chosen e.g. as $C = 2\sqrt{2}$ which does not depend on the order $k$ of the $U$-statistic $I_{n,k}(f)$.*

Formula (35) can be reformulated as $E(|Z_{\mu,k}(f)|)^{2M} \le E(\sigma \eta^k)^{2M}$, where $\eta$ is a standard normal random variable. Theorem 4.1 states that the tail distribution of $k!|Z_{\mu,k}(f)|$ satisfies an estimate similar to that of $\sigma|\eta|^k$. This can be deduced relatively simply from Proposition 5.1 and the Markov inequality $P(|Z_{\mu,k}(f)| > u) \le \frac{E(k!|Z_{\mu,k}(f)|)^{2M}}{u^{2M}}$ with an appropriate choice of the parameter $M$.

Proposition 5.2 gives a bound on the moments of $k! n^{-k/2} I_{n,k}(f)$ similar to the estimate (36) on the moments of $Z_{\mu,k}(f)$. The difference between them is that estimate (37) in Proposition 5.2 contains a factor $(1 + C\sqrt{\alpha})^{2kM}$ at its right-hand side, and it holds only for such moments $E \left( k! n^{-k/2} I_{n,k}(f) \right)^{2M}$ for which $k M_0 \le kM \le \alpha n \sigma^2$ with some constant $M_0$. The parameter $\alpha > 0$ in relation (36) can be chosen in an arbitrary way, but it yields a really useful estimate only for not too large values. Theorem 4.3 can be proved by means of the estimate in Proposition 5.2 and the Markov inequality. But because of the relatively weak estimate of Proposition 5.2 only the estimate of Theorem 4.3 can be proved for degenerate $U$-statistics. The main step both in the proof of Theorem 4.1 and 4.3 is to get good moment estimates.

A most important result of the probability theory, the so-called diagram formula about multiple Wiener–Itô integrals can be applied in the proof of Proposition 5.1. This result can be found e.g. in [17]. It enables us to rewrite the product of Wiener–Itô integrals as a sum of Wiener–Itô integrals of different order. It got the name 'diagram formula', because the kernel functions of the Wiener–Itô integrals appearing in the sum representation of the product of Wiener–Itô integrals are defined with the help of certain diagrams. As the expectation of a Wiener–Itô integral of order $k$ equals zero for all $k \ge 1$, the expectation of the product equals the sum of the constant terms (i.e. of the integrals of order zero) in the diagram formula. The sum of the constant terms in the diagram formula can be bounded, and such a calculation leads to the proof of Proposition 5.1.

A version of the diagram formula can be proved both for the product of multiple random integrals $J_{n,k}(f)$ defined in formula (2) (see [18]) or for degenerate $U$-statistics (see [20]) which expresses the product of multiple random integrals or degenerate $U$-statistics as a sum of multiple random integrals or degenerate $U$-statistics of different order. The main difference between these new and the original diagram formula about Wiener–Itô integrals is that in the case of random (non-Gaussian) integrals or degenerate $U$-statistics some new diagrams appear, and they give an additional contribution in the sum representation of the product of random integrals $J_{n,k}(f)$ or of degenerate $U$-statistics $I_{n,k}(f)$.

Proposition 5.2 can be proved by means of the diagram formula for the product of degenerate $U$-statistics and a good bound on the contribution of all integrals corresponding to the diagrams. Theorem 4.4 can be proved similarly by means of the diagram formula for the product of multiple random integrals $J_{n,k}(f)$ (see [18]). The main difficulty of such an approach arises, because the expected value of a $k$-fold random integral $J_{n,k}(f)$ (unlike that of a Wiener–Itô integral or degenerate $U$-statistic) may be non-zero also in the case $k \geq 1$. The expectation of all these integrals is small, but since the diagram formula contains a large number of such terms, it cannot supply such a sharp estimate for the moments random integrals $J_{n,k}(f)$ as we have for degenerate $U$-statistics $I_{n,k}(f)$. On the other hand, Theorem 4.4 can be deduced from Theorems 4.3, 3.2, and formulas (20) and (21).

*Remark:* The diagram formula is an important tool both in investigations in probability theory and statistical physics. The second chapter of the book [25] contains a detailed discussion of this formula. Paper [28] explains the combinatorial picture behind it, and it contains some interesting generalizations. Paper [31] is interesting because of a different reason. It shows how to prove central limit theorems for stationary processes in some non-trivial cases by means of the diagram formula. In this paper it is proved that the moments of the normalized partial sums have the right limit as the number of terms in them tends to infinity. Actually, the limit of the semi-invariants is investigated, but this can be considered as an appropriate reformulation of the study of the moments. The approach in paper [31] and the proof of the results mentioned in this work show some similarity, but there is also an essential difference between them. In paper [31] the limit of fixed moments is investigated, while e.g. in Problem a′) we want to get good asymptotics for such moments of $U$-statistics $I_{n,k}(f)$ whose order may depend on the sample size $n$ of the $U$-statistic. The reason behind this difference is that we want to get a good estimate of the probabilities defined in Problem a′) also for large numbers $u$, and this yields some large deviation character to the problem.

The statement of Example 4.2 follows relatively simply from another important result about multiple Wiener–Itô integrals, from the so-called Itô formula for multiple Wiener–Itô integrals (see e.g. [14] or [17]) which enables us to express the random integrals considered in Example 4.2 as the Hermite polynomial of an appropriately defined standard normal random variable.

Here I did not formulate the diagram formula, hence I cannot explain the details of the proof of Propositions 5.1 and 5.2. I discuss instead an analogous, but simpler problem briefly which may help in capturing the ideas behind the proofs outlined above.

Let us consider a sequence of independent and identically distributed random variables $\xi_1, \ldots, \xi_n$ with expectation zero, take their sum $S_n = \sum\limits_{j=1}^{n} \xi_j$, and let us try to give a good estimate on the moments $ES_n^{2M}$ for all $M = 1, 2, \ldots$. Because of the independence of the random variables $\xi_j$ and the condition $E\xi_j = 0$ we can write

$$ES_n^{2M} = \sum_{\substack{(j_1, \ldots, j_s, l_1, \ldots, l_s) \\ j_1 + \cdots + j_s = 2M, \, j_u \geq 2 \text{ for all } 1 \leq u \leq s, \\ l_u \neq l_{u'} \text{ if } u \neq u'}} E\xi_{l_1}^{j_1} \cdots E\xi_{l_s}^{j_s}. \quad (38)$$

Simple combinatorial considerations show that a dominating number of terms at the right-hand side of (38) are indexed by a vector $(j_1, \ldots, j_M, l_1, \ldots, l_M)$ such that $j_u = 2$ for all $1 \leq u \leq M$, and the number of such vectors is equal to $\binom{n}{M} \frac{(2M)!}{2^M} \sim n^M \frac{(2M)!}{2^M M!}$. The last asymptotic relation holds if the number $n$ of terms in the random sum $S_n$ is sufficiently large. The above considerations suggest that under not too restrictive conditions $ES_n^{2M} \sim (n\sigma^2)^M \frac{(2M)!}{2^M M!} = E\eta_{n\sigma^2}^{2M}$, where $\sigma^2 = E\xi^2$ is the variance of the terms in the sum $S_n$, and $\eta_u$ is a random variable with normal distribution with expectation zero and variance $u$. The question arises when the above heuristic argument gives a right estimate.

For the sake of simplicity let us restrict our attention to the case when the absolute value of the random variables $\xi_j$ is bounded by 1. Let us observe that even in this case we have to impose a condition that the variance $\sigma^2$ of the random variables $\xi_j$ is not too small. Indeed, let us consider such random variables $\xi_j$, for which $P(\xi_j = 1) = P(\xi_j = -1) = \frac{\sigma^2}{2}$, $P(\xi_j = 0) = 1 - \sigma^2$. These random variables $\xi_j$ have variance $\sigma^2$, and the contribution of the terms $E\xi_j^{2M}$, $1 \leq j \leq n$, to the sum in (38) equals $n\sigma^2$. If $\sigma^2$ is very small, then it may occur that $n\sigma^2 \gg (n\sigma^2)^M \frac{(2M)!}{2^M M!}$, and the approximation given for $ES_n^{2M}$ in the previous paragraph does not hold any longer. Let us observe that for larger moments $ES_n^{2M}$ the choice of a smaller variance $\sigma^2$ is sufficient to violate the asymptotic relation obtained by this approximation.

A similar picture arises in Proposition 5.2. If the variance of the random variable $I_{n,k}(f)$ is not too small, then those terms give the essential contribution to the moments of $I_{n,k}(f)$ which correspond to such diagrams which appear also in the diagram formula for Wiener–Itô integrals. The higher moment we estimate the stronger condition we have to impose on the variance of $I_{n,k}(f)$ to preserve this property and to get a good bound on the moment we consider.

In the next Section problems b), b′) and b″) will be discussed, where the distribution of the supremum of multiple random integrals $J_{n,k}(f)$, degenerate $U$-statistics $I_{n,k}(f)$ and multiple Wiener–Itô integrals $Z_{\mu,k}(f)$ will be estimated for an appropriate class of functions $f \in \mathcal{F}$. Under some appropriate conditions for the class of functions $\mathcal{F}$ a similar estimate can be proved in these problems as in their natural counterpart when only one function is taken. The only difference is that worse universal constants may appear in the new estimates. The conditions we had to impose in the results about problems a) and a′) appear in their counterparts problems b) and b') in a natural way. But these conditions also have some hidden, more surprising consequences in the study of the new problems.

## 6. On the supremum of random integrals and $U$-statistics

To formulate the results of this section first I introduce some notions which appear in their formulation. Such properties will be introduced which say about a class of functions that it has relatively small and in some sense dense finite subsets.

First I introduce the following definition.

**Definition of $L_p$-dense classes of functions with respect to some measure.** *Let us have a measurable space $(Y, \mathcal{Y})$ together with a $\sigma$-finite measure $\nu$ and a set $\mathcal{G}$ of $\mathcal{Y}$ measurable real valued functions on this space. For all $1 \le p < \infty$, we say that $\mathcal{G}$ is an $L_p$-dense class with respect to $\nu$ and with parameter $D$ and exponent $L$ if for all numbers $1 \ge \varepsilon > 0$ there exists a finite $\varepsilon$-dense subset $\mathcal{G}_\varepsilon = \{g_1, \ldots, g_m\} \subset \mathcal{G}$ in the space $L_p(Y, \mathcal{Y}, \nu)$ consisting of $m \le D\varepsilon^{-L}$ elements, i.e. there exists a set $\mathcal{G}_\varepsilon \subset \mathcal{G}$ with $m \le D\varepsilon^{-L}$ elements such that $\inf_{g_j \in \mathcal{G}_\varepsilon} \int |g - g_j|^p \, d\nu < \varepsilon^p$ for all functions $g \in \mathcal{G}$.*

The following notion will also be needed.

**Definition of $L_p$-dense classes of functions.** *Let us have a measurable space $(Y, \mathcal{Y})$ and a set $\mathcal{G}$ of $\mathcal{Y}$ measurable real valued functions on this space. We call $\mathcal{G}$ an $L_p$-dense class of functions, $1 \le p < \infty$, with parameter $D$ and exponent $L$ if it is $L_p$-dense with parameter $D$ and exponent $L$ with respect to all probability measures $\nu$ on $(Y, \mathcal{Y})$.*

The above introduced properties can be considered as possible versions of the so-called $\varepsilon$-entropy frequently applied in the literature. Nevertheless, there seems to exist no unanimously accepted version of this notion. Generally the above introduced definitions will be applied with the choice $p = 2$, but because of some arguments in this paper it was more natural to introduce them in a more general form. The first result I present can be considered as a solution of problem b'').

**Theorem 6.1.** *Let us consider a measurable space $(X, \mathcal{X})$ together with a $\sigma$-finite non-atomic measure $\mu$ on it, and let $\mu_W$ be a white noise with reference measure $\mu$ on $(X, \mathcal{X})$. Let $\mathcal{F}$ be a countable and $L_2$-dense class of functions $f(x_1, \ldots, x_k)$ on $(X^k, \mathcal{X}^k)$ with some parameter $D$ and exponent $L$ with respect to the product measure $\mu^k$ such that*

$$\int f^2(x_1, \ldots, x_k)\mu(\,dx_1)\ldots\mu(\,dx_k) \le \sigma^2 \quad \text{with some } 0 < \sigma \le 1 \text{ for all } f \in \mathcal{F}.$$

*Let us consider the multiple Wiener integrals $Z_{\mu,k}(f)$ introduced in formula (5) for all $f \in \mathcal{F}$. The inequality*

$$P\left(\sup_{f \in \mathcal{F}} |Z_{\mu,k}(f)| > u\right) \leq C(D+1)\exp\left\{-\alpha\left(\frac{u}{\sigma}\right)^{2/k}\right\} \quad \text{if } \left(\frac{u}{\sigma}\right)^{2/k} \geq ML\log\frac{2}{\sigma} \tag{39}$$

*holds with some universal constants $C = C(k) > 0$, $M = M(k) > 0$ and $\alpha = \alpha(k) > 0$.*

The next two results can be considered as a solution of problems b) and b′).

**Theorem 6.2.** *Let a probability measure $\mu$ be given on a measurable space $(X, \mathcal{X})$ together with a countable and $L_2$-dense class $\mathcal{F}$ of functions $f(x_1,\ldots,x_k)$ of $k$ variables with some parameter $D$ and exponent $L$, $L \geq 1$, on the product space $(X^k, \mathcal{X}^k)$ which satisfies the conditions*

$$\|f\|_\infty = \sup_{x_j \in X,\ 1 \leq j \leq k} |f(x_1,\ldots,x_k)| \leq 1, \qquad \text{for all } f \in \mathcal{F} \tag{40}$$

*and*

$$\|f\|_2^2 = Ef^2(\xi_1,\ldots,\xi_k)\ = \int f^2(x_1,\ldots,x_k)\mu(dx_1)\ldots\mu(dx_k) \leq \sigma^2 \tag{41}$$
$$\text{for all } f \in \mathcal{F}$$

*with some constant $0 < \sigma \leq 1$. Then there exist some constants $C = C(k) > 0$, $\alpha = \alpha(k) > 0$ and $M = M(k) > 0$ depending only on the parameter $k$ such that the supremum of the random integrals $J_{n,k}(f)$, $f \in \mathcal{F}$, defined by formula (2) satisfies the inequality*

$$P\left(\sup_{f \in \mathcal{F}} |J_{n,k}(f)| \geq u\right)\ \leq CD\exp\left\{-\alpha\left(\frac{u}{\sigma}\right)^{2/k}\right\} \tag{42}$$
$$\text{if}\quad n\sigma^2 \geq \left(\frac{u}{\sigma}\right)^{2/k} \geq M(L+\beta)^{3/2}\log\frac{2}{\sigma},$$

*where $\beta = \max\left(\frac{\log D}{\log n}, 0\right)$ and the numbers $D$ and $L$ agree with the parameter and exponent of the $L_2$-dense class $\mathcal{F}$.*

**Theorem 6.3.** *Let a probability measure $\mu$ be given on a measurable space $(X, \mathcal{X})$ together with a countable and $L_2$-dense class $\mathcal{F}$ of functions $f(x_1,\ldots,x_k)$ of $k$ variables with some parameter $D$ and exponent $L$, $L \geq 1$, on the product space $(X^k, \mathcal{X}^k)$ which satisfies conditions (40) and (41) with some constant $0 < \sigma \leq 1$. Beside these conditions let us also assume that the U-statistics $I_{n,k}(f)$ defined with the help of a sequence of independent $\mu$ distributed random variables $\xi_1,\ldots,\xi_n$ are degenerate for all $f \in \mathcal{F}$, or in an equivalent form, all functions $f \in \mathcal{F}$ are canonical with respect to the measure $\mu$. Then there exist some constants $C = C(k) > 0$, $\alpha = \alpha(k) > 0$ and $M = M(k) > 0$ depending* *only on the parameter $k$ such that the inequality*

$$P\left(\sup_{f\in\mathcal{F}} n^{-k/2}|I_{n,k}(f)|\geq u\right) \leq CD\exp\left\{-\alpha\left(\frac{u}{\sigma}\right)^{2/k}\right\} \tag{43}$$

$$if\quad n\sigma^2\geq\left(\frac{u}{\sigma}\right)^{2/k}\geq M(L+\beta)^{3/2}\log\frac{2}{\sigma},$$

*holds, where $\beta=\max\left(\frac{\log D}{\log n},0\right)$ and the number $D$ and $L$ agree with the parameter and exponent of the $L_2$-dense class $\mathcal{F}$.*

The above theorems whose proofs can be found in [19] or [22] in a more detailed version say that under some conditions on the class of functions $\mathcal{F}$ an almost as good estimate holds for problems b), b') and b'') as for the analogous problems a), a') and a''), where similar problems were investigated, but only one function $f$ was considered. An essential restriction in the results of Theorems 6.1, 6.2 and 6.3 is that the condition $\left(\frac{u}{\sigma}\right)^{2/k}\geq M(L,D)\log\frac{2}{\sigma}$ is imposed in them with some constant $M(L,D,k)$ depending on the exponent $L$ and parameter $D$ of the $L_2$-dense class $\mathcal{F}$. In Theorem 6.1 $M(L,D,k)=ML$ was chosen, in Theorems 6.2 and 6.3 $M(L,D,k)=M(L+\beta)^{3/2}$ with an appropriate universal constant $M=M(k)$ and $\beta=\max\left(0,\frac{\log D}{\log n}\right)$. We are interested not so much in a good choice of the quantity $M(L,D,k)$ in these results. Actually, they could have been chosen in a better way. We would like to understand why such conditions have to be imposed in these results.

I shall also discuss some other questions related to the above theorems. Beside the role of the lower bound on $\left(\frac{u}{\sigma}\right)^{2/k}$ one would also like to understand why we have imposed the condition of $L_2$-dense property for the class of functions $\mathcal{F}$ in Theorems 6.2 and 6.3. This is a stronger restriction than the condition about the $L_2$-dense property of the class $\mathcal{F}$ with respect to the measure $\mu^k$ imposed in Theorem 6.1. It may be a little bit mysterious why in Theorems 6.2 and 6.3 such a condition is needed by which this class of functions is $L_2(\nu)$-dense also with respect to such probability measures $\nu$ which seem to have no relation to our problems. I can give only a partial answer to this question. In the next section I present a very brief sketch of the proofs which shows that in the proof of Theorems 6.2 and 6.3 the $L_2$-dense property of the class of functions $\mathcal{F}$ is applied in the form as it was imposed. I shall discuss another question which also naturally arises in this context. One would like to know some results which enable us to check the $L_2$-dense property and show that it holds in many interesting cases.

I shall discuss still another problem related to the above results. One would like to weaken the condition by which the classes of functions $\mathcal{F}$ must be countable. Let me recall that in the Introduction I mentioned that our results can be applied in the study of some non-parametric maximum likelihood problems. In these applications such cases may occur where we have to work with the supremum of non-countably infinite random integrals. I shall discuss this question separately at the end of this section.

I show an example which shows that the condition $\left(\frac{u}{\sigma}\right)^{2/k} \geq M(L, D, k) \log \frac{2}{\sigma}$ with some appropriate constant $M(L, D, k) > 0$ cannot be omitted from Theorem 6.1. In this example $([0, 1], \mathcal{B})$, i.e. the interval $[0, 1]$ together with the Borel $\sigma$-algebra is taken as the measurable space $(X, \mathcal{X})$, and the Lebesgue measure $\lambda$ is considered on $[0, 1]$ together with the usual white noise $\lambda_W$ with the Lebesgue measure as its reference measure. Fix some number $\sigma > 0$, and define the class of functions of $k$ variables $\mathcal{F} = \mathcal{F}_\sigma$ on $([0, 1]^k, \mathcal{B}^k)$ as the indicator functions of the $k$-dimensional rectangles $\prod_{j=1}^{k} [a_j, b_j] \subset [0, 1]^k$ such that all numbers $a_j$ and $b_j$, $1 \leq j \leq k$, are rational, and the volume of these rectangles satisfy the condition $\prod_{j=1}^{k} (b_j - a_j) \leq \sigma^2$. It can be seen that this countable class of functions $\mathcal{F}$ is $L_2$-dense with respect to the measure $\lambda$, (moreover it is $L_2$-dense in the general sense), hence Theorem 6.1 can be applied to the supremum of the Wiener–Itô integrals $Z_{\lambda,k}(f)$ with the above class of functions $f \in \mathcal{F}$.

Let the above chosen number $\sigma > 0$ be sufficiently small and such that $\sigma^{2/k}$ is a rational number. Let us define $N = [\sigma^{-2/k}]$ functions $f_j \in \mathcal{F}$, where $[x]$ denotes the integer part of the number $x$ in the following way: The function $f_j$ is the indicator function of the $k$-dimensional cube we get by taking the $k$-fold direct product of the interval $[(j-1)\sigma^{2/k}, j\sigma^{2/k}]$ with itself, $1 \leq j \leq N$. Then all functions $f_j$ are elements of the above defined class of functions $\mathcal{F} = \mathcal{F}_\sigma$, and the Wiener–Itô integrals $Z_{\lambda,k}(f_j)$, $1 \leq j \leq N$, are independent random variables. Hence

$$P\left(\sup_{f \in \mathcal{F}} |Z_{\lambda,k}(f)| > u\right) \geq P\left(\sup_{1 \leq j \leq N} |Z_{\lambda,k}(f_j)| > u\right) = 1 - P(|Z_{\lambda,k}(f_1)| \leq u)^N \tag{44}$$

for all numbers $u > 0$. I will show with the help of relation (44) that for a small $\sigma > 0$ and such a number $u$ for which $\left(\frac{u}{\sigma}\right)^{2/k} = a \log \frac{2}{\sigma}$ with some $a < \frac{4}{k}$ the probability $P\left(\sup_{f \in \mathcal{F}} |Z_{\lambda,k}(f)| > u\right)$ is very close to 1.

By the Itô formula for multiple Wiener–Itô integrals (see e.g. [14]) the identity $Z_{\lambda,k}(f_j) = \sigma H_k(\eta_j)$ holds, where $H_k(\cdot)$ is the $k$-th Hermite polynomial with leading coefficient 1, and $\eta_j = \sigma^{-1/k} \int_{(j-1)\sigma^{2/k}}^{j\sigma^{2/k}} d\lambda_W$, hence it is a standard normal random variable. With the help of this relation it can be shown that for all $0 < \gamma < 1$ there exists some $\sigma_0 = \sigma_0(\gamma)$ such that $P(|Z_{\lambda,k}(f_1)| \leq u) \leq 1 - e^{-\gamma(u/\sigma)^{2/k}/2} = 1 - \left(\frac{\sigma}{2}\right)^{\gamma a/2}$ if $0 < \sigma < \sigma_0$. Hence relation (44) and the inequality $N \geq \sigma^{-2/k} - 1$ imply that $P\left(\sup_{f \in \mathcal{F}} |Z_{\lambda,k}(f)| > u\right) \geq 1 - \left(1 - \left(\frac{\sigma}{2}\right)^{\gamma a/2}\right)^{\sigma^{-2/k}-1}$. By choosing $\gamma$ sufficiently close to 1 it can be shown with the help of the above relation that with a sufficiently small $\sigma > 0$ and the above choice of the number $u$ the probability $P\left(\sup_{f \in \mathcal{F}} |Z_{\lambda,k}(f)| > u\right)$ is almost 1.

The above calculation shows that a condition of the type

$$\left(\frac{u}{\sigma}\right)^{2/k} \geq M(L, D, k) \log \frac{2}{\sigma}$$

is really needed in Theorem 6.1. With some extra work a similar example can be constructed in the case of Theorem 6.2. In this example the same space $(X, \mathcal{X})$ and the same class of functions $\mathcal{F} = \mathcal{F}_\sigma$ can be chosen, only the white noise has to be replaced for instance by a sequence of independent random variables $\xi_1, \ldots, \xi_n$ with uniform distribution on the unit interval and with a sufficiently large sample size $n$. (The lower bound on the sample size should depend also on $\sigma$.) Also in the case of Theorem 6.3 a similar example can be constructed. I omit the details.

The theory of Vapnik–Červonenkis classes is a fairly popular and important subject in probability theory. I shall show that this theory is also useful in the study of our problems. It provides a useful sufficient condition for the $L_2$-dense property of a class of functions, a property which played an important role in Theorems 6.2 and 6.3. To formulate the result interesting for us first I recall the notion of Vapnik–Červonenkis classes.

**Definition of Vapnik-Červonenkis classes of sets and functions.** *Let a set $S$ be given, and let us select a class $\mathcal{D}$ consisting of certain subsets of this set $S$. We call $\mathcal{D}$ a Vapnik–Červonenkis class if there exist two real numbers $B$ and $K$ such that for all positive integers $n$ and subsets $S_0(n) = \{x_1, \ldots, x_n\} \subset S$ of cardinality $n$ of the set $S$ the collection of sets of the form $S_0(n) \cap D$, $D \in \mathcal{D}$, contains no more than $Bn^K$ subsets of $S_0(n)$. We shall call $B$ the parameter and $K$ the exponent of this Vapnik–Červonenkis class.*

*A class of real valued functions $\mathcal{F}$ on a space $(Y, \mathcal{Y})$ is called a Vapnik–Červonenkis class if the collection of graphs of these functions is a Vapnik–Červonenkis class, i.e. if the sets $A(f) = \{(y, t) \colon y \in Y, \ \min(0, f(y)) \leq t \leq \max(0, f(y))\}$, $f \in \mathcal{F}$, constitute a Vapnik–Červonenkis class of subsets of the product space $S = Y \times R^1$.*

The theory about Vapnik–Červonenkis classes has generated a huge literature. Many sufficient conditions have been stated which ensure that certain classes of sets or functions are Vapnik–Červonenkis classes. Here I do not discuss them. I only present an important result of Richard Dudley, which states that a Vapnik–Červonenkis class of functions bounded by 1 is an $L_1$-dense class of functions.

**Theorem 6.4.** *Let $f(y)$, $f \in \mathcal{F}$, be a Vapnik–Červonenkis class of real valued functions on some measurable space $(Y, \mathcal{Y})$ such that $\sup\limits_{y \in Y} |f(y)| \leq 1$ for all $f \in \mathcal{F}$. Then $\mathcal{F}$ is an $L_1$-dense class of functions on $(Y, \mathcal{Y})$. More explicitly, if $\mathcal{F}$ is a Vapnik–Červonenkis class with parameter $B \geq 1$ and exponent $K > 0$, then it is an $L_1$-dense class with exponent $L = 2K$ and parameter $D = CB^2(4K)^{2K}$ with some universal constant $C > 0$.*

The proof of this result can be found in [28] (25 Approximation Lemma) or in my Lecture Note [22]. Formally, Theorem 6.4 gives a sufficient condition for a class of functions to be an $L_1$-dense class. But it is fairly simple to show that a class of functions satisfying the conditions of Theorem 6.4 is not only an $L_1$, but also an $L_2$-dense class. Indeed, an $L_1$-dense class of functions whose absolute values are bounded by 1 in the supremum norm is also an $L_2$-dense class, only with a possibly different exponent and parameter. I finish this section by discussing the problem how to replace the condition of countable cardinality of the class of functions in Theorems 6.2 and 6.3 by a useful weaker condition.

### 6.1. On the supremum of non-countable classes of random integrals and U-statistics

First I introduce the following notion.

**Definition of countably approximable classes of random variables.** *Let a class of random variables $U(f)$, $f \in \mathcal{F}$, indexed by a class of functions on a measurable space $(Y, \mathcal{Y})$ be given. We say that this class of random variables $U(f)$, $f \in \mathcal{F}$, is countably approximable if there is a countable subset $\mathcal{F}' \subset \mathcal{F}$ such that for all numbers $u > 0$ the sets $A(u) = \{\omega \colon \sup\limits_{f \in \mathcal{F}} |U(f)(\omega)| \geq u\}$ and $B(u) = \{\omega \colon \sup\limits_{f \in \mathcal{F}'} |U(f)(\omega)| \geq u\}$ satisfy the identity $P(A(u) \setminus B(u)) = 0$.*

It is fairly simple to see that in Theorems 6.1, 6.2 and 6.3 the condition about the countable cardinality of the class of functions $\mathcal{F}$ can be replaced by the weaker condition that the class of random variables $Z_{\mu,k}(f)$, $J_{n,k}(f)$ or $I_{n,k}(f)$, $f \in \mathcal{F}$, is a countably approximable class of functions. One would like to get some results which enable us to check this property. The following simple lemma (see Lemma 4.3 in [22]) may be useful for this.

**Lemma 6.5.** *Let a class of random variables $U(f)$, $f \in \mathcal{F}$, indexed by some set $\mathcal{F}$ of functions on a space $(Y, \mathcal{Y})$ be given. If there exists a countable subset $\mathcal{F}' \subset \mathcal{F}$ of the set $\mathcal{F}$ such that the sets $A(u) = \{\omega \colon \sup\limits_{f \in \mathcal{F}} |U(f)(\omega)| \geq u\}$ and $B(u) = \{\omega \colon \sup\limits_{f \in \mathcal{F}'} |U(f)(\omega)| \geq u\}$ introduced for all $u > 0$ in the definition of countable approximability satisfy the relation $A(u) \subset B(u - \varepsilon)$ for all $u > \varepsilon > 0$, then the class of random variables $U(f)$, $f \in \mathcal{F}$, is countably approximable.*

*The above property holds if for all $f \in \mathcal{F}$, $\varepsilon > 0$ and $\omega \in \Omega$ there exists a function $\bar{f} = \bar{f}(f, \varepsilon, \omega) \in \mathcal{F}'$ such that $|U(\bar{f})(\omega)| \geq |U(f)(\omega)| - \varepsilon$.*

Thus to prove the countable approximability property of a class of random variables $U(f)$, $f \in \mathcal{F}$, it is enough to check the condition formulated in the second paragraph of Lemma 6.5. I present an example when this condition can be checked. This example is particularly interesting, since in the study of non-parametric maximum likelihood problems such examples have to be considered.

Let us fix a function $f(x_1, \ldots, x_k)$, $\sup |f(x_1, \ldots, x_k)| \leq 1$, on the space $(X^k, \mathcal{X}^k) = (R^{ks}, \mathcal{B}^{ks})$ with some $s \geq 1$, where $\mathcal{B}^t$ denotes the Borel $\sigma$-algebra on the Euclidean space $R^t$, together with some probability measure $\mu$ on $(R^s, \mathcal{B}^s)$. For all vectors $(u_1, \ldots, u_k)$, $(v_1, \ldots, v_k)$ such that $u_j, v_j \in R^s$ and $u_j \leq v_j$, $1 \leq j \leq k$, (i.e. all coordinates of $u_j$ are smaller than or equal to the corresponding coordinate of $v_j$) let us define the function $f_{u_1, \ldots, u_k, v_1, \ldots, v_k}$ which equals the function $f$ on the rectangle $[u_1, v_1] \times \cdots [u_k, v_k]$, and it is zero outside of this rectangle.

Let us consider a sequence of i.i.d. random variables $\xi_1, \ldots, \xi_n$ taking values in the space $(R^s, \mathcal{B}^s)$ with some distribution $\mu$, and define the empirical measure $\mu_n$ and random integrals $J_{n,k}(f_{u_1, \ldots, u_k, v_1, \ldots, v_k})$ by formulas (1) and (2) for all vectors $(u_1, \ldots, u_k)$, $(v_1, \ldots, v_k)$, $u_j \leq v_j$ for all $1 \leq j \leq k$, with the above defined functions $f_{u_1, \ldots, u_k, v_1, \ldots, v_k}$. The following result holds (see Lemma 4.4 in [22]).

**Lemma 6.6.** *Let us take $n$ independent and identically distributed random variables $\xi_1, \ldots, \xi_n$ with values in the space $(R^s, \mathcal{B}^s)$. Let us define with the help of their distribution $\mu$ and the empirical distribution $\mu_n$ determined by them the class of random variables $J_{n,k}(f_{u_1, \ldots, u_k, v_1, \ldots, v_k})$ introduced in formula (2), where the class of kernel functions $\mathcal{F}$ in these integrals consists of all functions $f_{u_1, \ldots, u_k, v_1, \ldots, v_k} \in (R^{sk}, \mathcal{B}^{sk})$, $u_j, v_j \in R^s$, $u_j \leq v_j$, $1 \leq j \leq k$, introduced in the last but one paragraph. This class of random variables $J_{n,k}(f)$, $f \in \mathcal{F}$, is countably approximable.*

Let me also remark that the class of functions $f_{u_1, \ldots, u_k, v_1, \ldots, v_k}$ is also an $L_2$-dense class of functions, actually it is also a Vapnik–Červonenkis class of functions. As a consequence, Theorem 6.2 can be applied to this class of functions.

To clarify the background of the above results I make the following remark. The class of random variables $Z_{\mu,k}(f)$, $J_{n,k}(f)$ or $I_{n,k}(f)$, $f \in \mathcal{F}$, can be considered as a stochastic process indexed by the functions $f \in \mathcal{F}$, and we estimate the supremum of this stochastic process. In the study of a stochastic process with a large parameter set one introduces some smoothness type property of the trajectories which can be satisfied. Here we followed a very similar approach. The condition formulated in the second paragraph of Lemma 6.5 can be considered as the smoothness type property needed in our problem.

In the study of a general stochastic process one has to make special efforts to find its right version with sufficiently smooth trajectories. In the case of the random processes $J_{n,k}(f)$ or $I_{n,k}(f)$, $f \in \mathcal{F}$, this right version can be constructed in a natural, simple way. A finite sequence of random variables $\xi_1(\omega), \ldots, \xi_n(\omega)$ is given at the start, and the random integrals $J_{n,k}(f)(\omega)$ or $U$-statistics $I_{n,k}(f)(\omega)$, $f \in \mathcal{F}$, can be constructed separately for all $\omega \in \Omega$ on the probability field $(\Omega, \mathcal{A}, P)$ where the random variables $\xi_1(\omega), \ldots, \xi_n(\omega)$ are living. It has to be checked whether the 'trajectories' of this random process have the 'smoothness properties' necessary for us. The case of a class of Wiener–Itô integrals $Z_{\mu,k}(f)$, $f \in \mathcal{F}$, is different. Wiener–Itô integrals are defined with the help of some $L_2$-limit procedure. Hence each random integral $Z_{\mu,k}(f)$ is defined only with probability 1, and in the case of a non-countable set of functions $\mathcal{F}$ the right version $Z_{\mu,k}(f)$, $f \in \mathcal{F}$, of the Wiener–Itô integrals has to be found to get a countably approximable class of random variables.

R. M. Dudley (see e.g. [5]) worked out a rather deep theory to overcome the measurability difficulties appearing in the case of a non-countable set of random variables by working with analytic sets, Suslin property, outer probability, and so on. I must admit that I do not know the precise relation between this theory and our method. At any rate, in the problems discussed here our elementary approach seems to be satisfactory.

In the next two sections I discuss the idea of the proof of Theorems 6.1, 6.2 and 6.3. A simple and natural approach, the so-called chaining argument suffices to prove Theorem 6.1. In the case of Theorems 6.2 and 6.3 this chaining argument can only help to reduce the proof to a slightly weaker statement, and we apply an essentially different method based on some randomization arguments to complete the proof. Since in the multivariate case $k \geq 2$ some essential additional difficulties appear, it seemed to be more natural to discuss it in a separate section.

## 7. The method of proof of Theorems 6.1, 6.2 and 6.3

There is a simple but useful method, called the chaining argument, which helps to prove Theorem 6.1. It suggests to take an appropriate increasing sequence $\mathcal{F}_j$, $j = 0, 1, \ldots$, of $L_2$-dense subsets of the class of functions $\mathcal{F}$ and to estimate the supremum of the Wiener–Itô integrals $Z_{\mu,k}(f)$, $f \in \mathcal{F}_j$, for all $j = 0, 1, \ldots$.

In the application of this method first we define a sequence of subclasses $\mathcal{G}_j$ of $\mathcal{F}$, $j = 0, 1, 2, \ldots$, such that $\mathcal{G}_j = \{g_{j,1}, \ldots, g_{j,m_j}\} \subset \mathcal{F}$ is an $2^{-jk}\sigma$-dense subset of $\mathcal{F}$ in the $L_2(\mu^k)$-norm, i.e. they satisfy the relation

$$
\begin{aligned}
&\inf_{1 \leq l \leq m_j} \rho(f, g_{j,l})^2 \\
&= \inf_{1 \leq l \leq m_j} \int (f(x_1, \ldots, x_k) - g_{j,l}(x_1, \ldots, x_k))^2 \mu(dx_1) \ldots \mu(dx_k) \\
&\leq 2^{-2jk}\sigma^2
\end{aligned}
\tag{45}
$$

for all $f \in \mathcal{F}$, and also the inequality $m_j \leq D2^{jkL}\sigma^{-L}$ holds. Such sets $\mathcal{G}_j$ exist because of the conditions of Theorem 6.1. Let us also define the classes of functions $\mathcal{F}_j = \bigcup_{p=0}^{j} \mathcal{G}_p$, and sets

$$
B_j = B_j(u) = \left\{ \omega \colon \sup_{f \in \mathcal{F}_j} |Z_{\mu,k}(f)(\omega)| \geq u\left(1 - 2^{-jk/2}\right) \right\}, \quad j = 0, 1, 2, \ldots.
$$

Given a function $f_{j+1,l} \in \mathcal{G}_{j+1}$ let us choose such a function $f_{j,l'} \in \mathcal{F}_j$ with some $l' = l'(l)$ for which $\rho(f_{j,l'}, f_{j+1,l}) \leq 2^{-jk}\sigma$ with the function $\rho(f, g)$ defined in

formula (45). Then

$$P(B_{j+1}) \leq P(B_j) + \sum_{l=1}^{m_{j+1}} P\left(|Z_{\mu,k}(f_{j+1,l} - f_{j,l'})| > u2^{-k(j+1)/2}\right). \quad (46)$$

Theorem 4.1 yields a good estimate of the terms in the sum at the right-hand side of (46), and it also provides a good bound of the probability $P(B_0)$. With the help of some small modification of the construction it can be achieved that also the relation $\bigcup_{j=0}^{\infty} \mathcal{F}_j = \mathcal{F}$ holds. The proof of Theorem 6.1 follows from the estimates obtained in such a way.

Theorem 6.2 can be deduced from Theorem 6.3 relatively simply with the help of Theorem 3.3, since Theorem 6.3 enables us to give a good bound on all terms in the sum at the right-hand side of formula (19). The only non-trivial step in this argument is to show that the set of functions $f_V$, $f \in \mathcal{F}$, appearing in formula (19) satisfy the estimates needed in the application of Theorem 6.3. Relations (20) and (21) are parts of the needed estimates. Beside this, it has to be shown that if $\mathcal{F}$ is an $L_2$-dense class of functions, then the same relation holds for the classes of functions $\mathcal{F}_V = \{f_V : f \in \mathcal{F}\}$ for all sets $V \subset \{1, \ldots, k\}$. This relation can also be shown with the help of a not too difficult proof (see [19] or [22]), but this question will be not discussed here.

One may try to prove Theorem 6.3, similarly to Theorem 6.1, with the help of the chaining argument. But this method does not work well in this case. The reason for its weakness is that the tail distribution of a degenerate $U$-statistic with a small variance $\sigma^2$ does not satisfy such a good estimate as the tail distribution of a multiple Wiener–Itô integral. At this point the condition $u \leq n^{k/2}\sigma^{k+1}$ in Theorem 4.2 plays an important role. Let us recall that, as Example 4.5 shows, the tail distribution of the normalized degenerated $U$-statistics $n^{-k/2}I_{n,k}(f)$ satisfies only a relatively weak estimate at level $u$ if $u \gg n^{k/2}\sigma^{k+1}$. We may try to work with an estimate analogous to relation (46) in the proof of Theorem 6.3. But the probabilities appearing at the right-hand side of such an estimate cannot be well estimated for large indices $j$.

Thus we can start the procedure of the chaining argument, but after finitely many steps we have to stop it. In such a way we can find a relatively dense subset $\mathcal{F}_0 \subset \mathcal{F}$ (in $L_2(\mu)$ norm) such that a good estimate can be given for the distribution of the supremum $\sup_{f \in \mathcal{F}_0} I_{n,k}(f)$. This result enables us to reduce Theorem 6.3 to a slightly weaker statement formulated in Proposition 7.1 below, but it yields no more help. Nevertheless, such a reduction is useful.

**Proposition 7.1.** *Let us have a probability measure $\mu$ on a measurable space $(X, \mathcal{X})$ together with a sequence of independent and $\mu$ distributed random variables $\xi_1, \ldots, \xi_n$ and a countable $L_2$-dense class $\mathcal{F}$ of canonical kernel functions $f = f(x_1, \ldots, x_k)$ (with respect to the measure $\mu$) with some parameter $D$ and exponent $L$ on the product space $(X^k, \mathcal{X}^k)$ such that all functions $f \in \mathcal{F}$ satisfy conditions (40) and (41) with some $0 < \sigma \leq 1$. Let us consider the (degenerate) $U$-statistics $I_{n,k}(f)$ with the random sequence $\xi_1, \ldots, \xi_n$ and kernel functions* $f \in \mathcal{F}$. *There exists a sufficiently large constant $K = K(k)$ together with some numbers $\bar{C} = \bar{C}(k) > 0$, $\gamma = \gamma(k) > 0$ and threshold index $A_0 = A_0(k) > 0$ depending only on the order $k$ of the U-statistics such that if $n\sigma^2 > K(L+\beta)\log n$ with $\beta = \max\left(\frac{\log D}{\log n}, 0\right)$, then the degenerate U-statistics $I_{n,k}(f)$, $f \in \mathcal{F}$, satisfy the inequality*

$$P\left(\sup_{f \in \mathcal{F}} |n^{-k/2}I_{n,k}(f)| \geq An^{k/2}\sigma^{k+1}\right) \leq \bar{C}e^{-\gamma A^{1/2k}n\sigma^2} \quad \text{if } A \geq A_0. \quad (47)$$

The statement of Proposition 7.1 is similar to that of Theorem 6.3. The essential difference between them is that Proposition 7.1 yields an estimate only for $u \geq A_0 n^{k/2}\sigma^{k+1}$ with a sufficiently large constant $A_0$, i.e. for relatively large numbers $u$. In the case $u \gg n^{k/2}\sigma^{k+1}$ it yields a weaker estimate than formula (43) in Theorem 6.3, but actually we need this estimate only in the case of the number $A$ in formula (47) being bounded away both from zero and infinity.

The proof of Proposition 7.1, briefly explained below, is based on an inductive procedure carried out by means of a symmetrization argument. In each step of this induction we diminish the number $A_0$ for which we show that inequality (47) holds for all numbers $An^{k/2}\sigma^{k+1}$ with $A \geq A_0$. This diminishing of the number $A_0$ is done as long as it is possible. It has to be stopped at such a number $A_0$ for which the probability $P(|n^{-k/2}I_{n,k}(f)| \geq A_0 n^{k/2}\sigma^{k+1})$ can be well estimated by Theorem 4.3 for all functions $f \in \mathcal{F}$. This has the consequence that Proposition 7.1 yields just such a strong estimate which is needed to reduce the proof of Theorem 6.3 to a statement that can be proved by means of the chaining argument.

In the symmetrization argument applied in the proof of Proposition 7.1 several additional difficulties arise if the multivariate case $k \geq 2$ is considered. Hence in this section only the case $k = 1$ is discussed. A degenerate U-statistic $I_{n,1}(f)$ of order 1 is the sum of independent, identically distributed random variables with expectation zero. In this paper the proof of Proposition 7.1 will be only briefly explained. A detailed proof can be found in [19] or [22]. Let me also remark that the method of these works was taken from Alexander's paper [2], where all ideas appeared in a different context.

We shall bound the probability appearing at the left-hand side of (47) (if $k = 1$) from above by the probability of the event that the supremum of appropriate randomized sums is larger than some number. We apply a symmetrization method which means that we estimate the expression we want to bound by means of a randomized (symmetrized) expression. Lemma 7.2, formulated below, has such a character.

**Lemma 7.2.** *Let a countable class of functions $\mathcal{F}$ on a measurable space $(X, \mathcal{X})$ and a real number $0 < \sigma < 1$ be given. Consider a sequence of independent, identically distributed $X$-valued random variables $\xi_1, \ldots, \xi_n$ such that $Ef(\xi_1) = 0$, $Ef^2(\xi_1) \leq \sigma^2$ for all $f \in \mathcal{F}$ together with another sequence $\varepsilon_1, \ldots, \varepsilon_n$ of independent random variables with distribution $P(\varepsilon_j = 1) = P(\varepsilon_j = -1) = \frac{1}{2}$,*

$1 \leq j \leq n$, *independent also of the random sequence* $\xi_1, \ldots, \xi_n$. *Then*

$$P \left( \frac{1}{\sqrt{n}} \sup_{f \in \mathcal{F}} \left| \sum_{j=1}^{n} f(\xi_j) \right| \geq A n^{1/2} \sigma^2 \right)$$

$$\leq 4P \left( \frac{1}{\sqrt{n}} \sup_{f \in \mathcal{F}} \left| \sum_{j=1}^{n} \varepsilon_j f(\xi_j) \right| \geq \frac{A}{3} n^{1/2} \sigma^2 \right) \quad \text{if } A \geq \frac{3\sqrt{2}}{\sqrt{n}\sigma}. \tag{48}$$

Let us first understand why Lemma 7.2 can help in the proof of Proposition 7.1. It enables to reduce the estimate of the probability at the left-hand side of formula (48) to that at its right-hand side. This reduction turned out to be useful for the following reason. At the right-hand side of formula 7.4 the probability of such an event appears which depends on the random variables $\xi_1, \ldots, \xi_n$ and some randomizing terms $\varepsilon_1, \ldots, \varepsilon_n$. Let us estimate the probability of this event by bounding first its conditional probability under the condition that the values of the random variables $\xi_1, \ldots, \xi_n$ are prescribed. These conditional probabilities can be well estimated by means of Hoeffding's inequality formulated below, and the estimates we get for them also yield a good bound on the expression at the right-hand side of (48).

Hoeffding's inequality, (see e.g. in [28] pp. 191–192), more precisely its special case we need here, states that the linear combinations of independent random variables $\varepsilon_j$, $P(\varepsilon_j = 1) = P(\varepsilon_j = -1) = \frac{1}{2}$, $1 \leq j \leq n$, behave so as the central limit theorem suggests. More explicitly, the following inequality holds.

**Theorem 7.3 (Hoeffding's inequality).** *Let* $\varepsilon_1, \ldots, \varepsilon_n$ *be independent random variables,* $P(\varepsilon_j = 1) = P(\varepsilon_j = -1) = \frac{1}{2}$, $1 \leq j \leq n$, *and let* $a_1, \ldots, a_n$ *be arbitrary real numbers. Put* $V = \sum_{j=1}^{n} a_j \varepsilon_j$. *Then*

$$P(V > y) \leq \exp \left\{ -\frac{y^2}{2 \sum_{j=1}^{n} a_j^2} \right\} \quad \text{for all } y > 0. \tag{49}$$

As we shall see, the application of Lemma 7.2 together with the above mentioned conditioning argument and Hoeffding's inequality enable us to reduce the estimation of the distribution of $\sup_{f \in \mathcal{F}} \sum_{j=1}^{n} f(\xi_j)$ to that of $\sup_{f \in \mathcal{F}} \sum_{j=1}^{n} f^2(\xi_j) = \sup_{f \in \mathcal{F}} \sum_{j=1}^{n} [f^2(\xi_j) - Ef^2(\xi_j)] + n \sup_{f \in \mathcal{F}} Ef^2(\xi_1)$. At first sight it may seem so that we did not gain very much by applying this approach. The estimation of the supremum of a class of sums of independent and identically distributed random variables was replaced by the estimation of a similar supremum. But a closer look shows that this method can help us in finding a proof of Proposition 7.1. We have to follow at what level we wanted to bound the distribution of the supremum in the original problem, and what level we have to choose in the modified problem to get a good estimate in the problem we are interested in. It turns out that in the second problem we need a good estimate about the distribution of the supremum of a class of sums of independent and identically distributed random variables at a considerable higher level. This observation enables us to work out an inductive procedure which leads to the proof of Proposition 7.1.

Indeed, in Proposition 7.1 estimate (47) has to be proved for all numbers $A \geq A_0$ with some appropriate number $A_0$. This estimate trivially holds if $A_0 > \sigma^{-2}$, because in this case condition (40) about the functions $f \in \mathcal{F}$ implies that the probability at the left-hand side of (47) equals zero. The argument of the previous paragraph suggests the following statement: If relation (47) holds for some constant $A_0$, then it also holds for a smaller $A_0$. Hence Proposition 7.1 can be proved by means of an inductive procedure in which the number $A_0$ is diminished at each step.

The actual proof consists of an elaboration of the details in the above heuristic approach. An inductive procedure is applied in which it is shown that if relation (47) holds with some number $A_0$ for a class of functions $\mathcal{F}$ satisfying the conditions of Proposition 7.1, then this relation also holds for it if $A_0$ is replaced by $A_0^{3/4}$, provided that $A_0$ is larger than some fixed universal constant. I would like to emphasize that we prove this statement not only for the class of functions $\mathcal{F}$ we are interested in, but simultaneously for all classes of functions which satisfy the conditions of Proposition 7.1. As we want to prove the inductive statement for a class of functions $\mathcal{F}$, then we apply our previous information not to this class, but to another appropriately defined class of functions $\mathcal{F}' = \mathcal{F}'(\mathcal{F})$ which also satisfies the conditions of Propositions 7.1. I omit the details of the proof, I only discuss one point which deserves special attention.

Hoeffding's inequality, applied in the justification of the inductive procedure leading to the proof of Proposition 7.1 gives an estimate for the distribution of a single sum, while we need a good estimate on the supremum of a class of sums. The question may arise whether this does not cause some problem in the proof. I try to briefly explain that the reason to introduce the condition about the $L_2$-dense property of the class $\mathcal{F}$ was to overcome this difficulty.

In the inductive procedure we want to prove that relation (47) holds for all $A \geq A_0^{3/4}$ if it holds for all $A \geq A_0$. It can be shown by means of the inductive assumption which states that relation (47) holds for $A \geq A_0$ and Hoeffding's inequality Theorem 7.3 that there is a set $D \subset \Omega$ such that the conditional probabilities

$$P\left(\frac{1}{\sqrt{n}}\left|\sum_{j=1}^{n}\varepsilon_j f(\xi_j)\right| \geq \frac{An^{1/2}\sigma^2}{6}\,\middle|\,\xi_1(\omega),\ldots\xi_n(\omega)\right) \tag{50}$$

are very small for all $f \in \mathcal{F}$, and the probability of the set $\Omega \setminus D$ is negligibly small. Let me emphasize that at this step of the proof we can give a good estimate about the conditional probability in (50) for all functions $f \in \mathcal{F}$ if $\omega \in D$, but we cannot work with their supremum which we would need to apply formula (48). This difficulty can be overcome with the help of the following argument.

Let us introduce the (random) probability measure $\nu = \nu(\omega)$ uniformly distributed on the points $\xi_1(\omega), \ldots, \xi_n(\omega)$ for all $\omega \in D$. Let us observe that the (random) measure $\nu$ has a support consisting of $n$ points, and the $\nu$-measure of all points in the support of $\nu$ equals $\frac{1}{n}$. This implies that the supremum of a function defined on the support of the measure $\nu$ can be bounded by means of the $L_2(\nu)$-norm of this function. This property together with the $L_2(\nu)$-dense property of the class of functions $\mathcal{F}$ imposed in the conditions of Proposition 7.1 imply that a finite set $\{f_1, \ldots, f_m\} \subset \mathcal{F}$ can be chosen with relatively few elements $m$ in such a way that for all $f \in \mathcal{F}$ there is some function $f_l$, $1 \le l \le m$, whose distance from the function $f$ in the $L_2(\nu)$ norm is less than $A\sigma^2/6$, hence $\inf_{1 \le l \le m} n^{-1/2} \left| \sum_{j=1}^{n} \varepsilon_j (f(\xi_j) - f_l(\xi_j)) \right| \le n^{1/2} \int |f - f_l| d\nu \le \frac{An^{1/2}\sigma^2}{6}$. The condition that $\mathcal{F}$ is $L_2$-dense with exponent $L$ and parameter $D$ enables us to give a good upper bound on the number $m$. This is the point, where the condition that the class of functions $\mathcal{F}$ is $L_2$-dense was exploited in its full strength. Since we can give a good bound on the conditional probability in (50) for all functions $f = f_l$, $1 \le l \le m$, we can bound the probability at the right-hand side of (48). It turns out that the estimate we get in such a way is sufficiently sharp, and the inductive statement, hence also Proposition 7.1 can be proved by working out the details.

I briefly explain the proof of Lemma 7.2. The randomizing terms $\varepsilon_j$, $1 \le j \le n$, in it can be introduced with the help of the following simple lemma.

**Lemma 7.4.** *Let $\xi_1, \ldots, \xi_n$ and $\bar\xi_1, \ldots, \bar\xi_n$ be two sequences of independent and identically distributed random variables with the same distribution $\mu$ on some measurable space $(X, \mathcal{X})$, independent of each other. Let $\varepsilon_1, \ldots, \varepsilon_n$ be a sequence of independent random variables $P(\varepsilon_j = 1) = P(\varepsilon_j = -1) = \frac{1}{2}$, $1 \le j \le n$, which is independent of the random sequences $\xi_1, \ldots, \xi_n$ and $\bar\xi_1, \ldots, \bar\xi_n$. Take a countable set of functions $\mathcal{F}$ on the space $(X, \mathcal{X})$. Then the set of random variables*

$$\frac{1}{\sqrt{n}} \sum_{j=1}^{n} \left( f(\xi_j) - f(\bar\xi_j) \right), \quad f \in \mathcal{F},$$

*and its randomized version*

$$\frac{1}{\sqrt{n}} \sum_{j=1}^{n} \varepsilon_j \left( f(\xi_j) - f(\bar\xi_j) \right), \quad f \in \mathcal{F},$$

*have the same joint distribution.*

Lemma 7.2 can be proved by means of Lemma 7.4 and some calculations. There is one harder step in the calculations. A probability of the type

$$P \left( \frac{1}{\sqrt{n}} \sup_{f \in \mathcal{F}} \sum_{j=1}^{n} f(\xi_j) > u \right)$$

has to be bounded from above by means of a probability of the type

$$P\left(\frac{1}{\sqrt{n}}\sup_{f\in\mathcal{F}}\sum_{j=1}^{n}\left(f(\xi_j)-f(\bar{\xi}_j)\right)>u-K\right)$$

with some number $K>0$. (Here the notation of Lemma 7.4 is applied.) At this point the following symmetrization lemma may be useful.

**Lemma 7.5 (Symmetrization Lemma).** *Let $Z_p$ and $\bar{Z}_p$, $p=1,2,\dots$, be two sequences of random variables independent of each other, and let the random variables $\bar{Z}_p$, $p=1,2,\dots$, satisfy the inequality*

$$P(|\bar{Z}_p|\le\alpha)\ge\beta \quad \text{for all } p=1,2,\dots \tag{51}$$

*with some numbers $\alpha\ge 0$ and $\beta\ge 0$. Then*

$$P\left(\sup_{1\le p<\infty}|Z_p|>\alpha+u\right)\le\frac{1}{\beta}P\left(\sup_{1\le p<\infty}|Z_p-\bar{Z}_p|>u\right) \quad \text{for all } u>0.$$

The proof of Lemma 7.5 can be found for instance in [28] (8 Symmetrization Lemma) or in [22] Lemma 7.1.

Let us list the element of the countable class of functions $\mathcal{F}$ in Lemma 7.2 in the form $\mathcal{F}=\{f_1,f_2,\dots,\}$. Then Lemma 7.2 can be proved by means of Lemmas 7.4 and 7.5 with the choice of the random variables

$$Z_p=\frac{1}{\sqrt{n}}\sum_{j=1}^{n}f_p(\xi_j) \quad \text{and} \quad \bar{Z}_p=\frac{1}{\sqrt{n}}\sum_{j=1}^{n}f_p(\bar{\xi}_j), \quad p=1,2,\dots. \tag{52}$$

I omit the details.

One may try to generalize the above sketched proof of Theorem 6.3 to the multivariate case $k\ge 2$. Here the question arises on how to generalize Lemma 7.2 to the multivariate case and how to prove this generalization. These are highly non-trivial problems. This will be the main subject of the next section.

## 8. On the proof of Theorem 6.3 in the multivariate case

Here we are mainly interested in the question how to carry out the symmetrization procedure in the proof of Proposition 7.1 to the multivariate case $k\ge 2$. It turned out that it is possible to reduce this problem to the investigation of modified $U$-statistics, where $k$ independent copies of the original random sequence are taken and put into the $k$ different arguments of the kernel function of the $U$-statistic of order $k$. Such modified versions of $U$-statistics are called decoupled $U$-statistics in the literature, and they can be better studied by means of the symmetrization argument we are going to apply. To give a precise meaning of the above statements some definitions have to be introduced and some results have to be formulated. I introduce the following notions.

**The definition of decoupled and randomized decoupled $U$-statistics.**
*Let us have $k$ independent copies $\xi_1^{(j)}, \ldots, \xi_n^{(j)}$, $1 \leq j \leq k$, of a sequence $\xi_1, \ldots, \xi_n$ of independent and identically distributed random variables taking their values in a measurable space $(X, \mathcal{X})$ together with a measurable function $f(x_1, \ldots, x_k)$ on the product space $(X^k, \mathcal{X}^k)$ with values in a separable Banach space. Then the decoupled $U$-statistic determined by the random sequences $\xi_1^{(j)}, \ldots, \xi_n^{(j)}$, $1 \leq j \leq k$, and kernel function $f$ is defined by the formula*

$$\bar{I}_{n,k}(f) = \frac{1}{k!} \sum_{\substack{1 \leq l_j \leq n, \, j = 1, \ldots, k \\ l_j \neq l_{j'} \text{ if } j \neq j'}} f\left(\xi_{l_1}^{(1)}, \ldots, \xi_{l_k}^{(k)}\right). \tag{53}$$

*Let us have beside the sequences $\xi_1^{(j)}, \ldots, \xi_n^{(j)}$, $1 \leq j \leq k$, and function $f(x_1, \ldots, x_k)$ a sequence of independent random variables $\varepsilon = (\varepsilon_1, \ldots, \varepsilon_n)$, $P(\varepsilon_l = 1) = P(\varepsilon_l = -1) = \frac{1}{2}$, $1 \leq l \leq n$, which is independent also of the sequences of random variables $\xi_1^{(j)}, \ldots, \xi_n^{(j)}$, $1 \leq j \leq k$. We define the randomized decoupled $U$-statistic determined by the random sequences $\xi_1^{(j)}, \ldots, \xi_n^{(j)}$, $1 \leq j \leq k$, the kernel function $f$ and the randomizing sequence $\varepsilon_1, \ldots, \varepsilon_n$ by the formula*

$$\bar{I}_{n,k}^{\varepsilon}(f) = \frac{1}{k!} \sum_{\substack{1 \leq l_j \leq n, \, j = 1, \ldots, k \\ l_j \neq l_{j'} \text{ if } j \neq j'}} \varepsilon_{l_1} \cdots \varepsilon_{l_k} f\left(\xi_{l_1}^{(1)}, \ldots, \xi_{l_k}^{(k)}\right). \tag{54}$$

Our first goal is to reduce the study of inequality (47) in Proposition 7.1 to an analogous problem about the supremum of decoupled $U$-statistics defined above. Then we want to show that a symmetrization argument enables us to reduce this problem to the study of randomized decoupled $U$-statistics introduced in formula (54). A result of de la Peña and Montgomery–Smith formulated below helps to carry out such a program. Let me remark that both in the definition of decoupled $U$-statistics and in the result of de la Peña and Montgomery–Smith functions $f$ taking their values in a separable Banach space were considered, i.e. we did not restrict our attention to real-valued functions. This choice was motivated by the fact that in such a general setting we can get a simpler proof of inequality (56) presented below. (The definition of $U$-statistics given in formula (3) is also meaningful in the case of Banach-space valued functions $f$.)

**Theorem 8.1 (de la Peña and Montgomery–Smith).** *Let us consider a sequence of independent and identically distributed random variables $\xi_1, \ldots, \xi_n$ on a measurable space $(X, \mathcal{X})$ together with $k$ independent copies $\xi_1^{(j)}, \ldots, \xi_n^{(j)}$, $1 \leq j \leq k$. Let us also have a function $f(x_1, \ldots, x_k)$ on the $k$-fold product space $(X^k, \mathcal{X}^k)$ which takes its values in a separable Banach space $B$. Define the $U$-statistic and decoupled $U$-statistic $I_{n,k}(f)$ and $\bar{I}_{n,k}(f)$ with the help of the above random sequences $\xi_1, \ldots, \xi_n$, $\xi_1^{(j)}, \ldots, \xi_n^{(j)}$, $1 \leq j \leq k$, and kernel function $f$.*

There exist some constants $\bar{C} = \bar{C}(k) > 0$ and $\gamma = \gamma(k) > 0$ depending only on the order $k$ of the $U$-statistic such that

$$P\left(\|I_{n,k}(f)\| > u\right) \leq \bar{C}P\left(\|\bar{I}_{n,k}(f)\| > \gamma u\right) \tag{55}$$

for all $u > 0$. Here $\|\cdot\|$ denotes the norm in the Banach space $B$ where the function $f$ takes its values.

More generally, if we have a countable sequence of functions $f_s$, $s = 1, 2, \ldots$, taking their values in the same separable Banach-space, then

$$P\left(\sup_{1 \leq s < \infty} \|I_{n,k}(f_s)\| > u\right) \leq \bar{C}P\left(\sup_{1 \leq s < \infty} \|\bar{I}_{n,k}(f_s)\| > \gamma u\right). \tag{56}$$

The proof of Theorem 8.1 can be found in [4] or in Appendix B of my Lecture Note [22]. Actually [4] contains only the proof of inequality (55), but (56) can be deduced from it simply by introducing appropriate separable Banach spaces and by exploiting that the universal constants in formula (55) do not depend on the Banach space where the random variables are living. Theorem 8.1 is useful for us, because it shows that Proposition 7.1 simply follows from its version presented in Proposition 8.2 below, where $U$-statistics are replaced by decoupled $U$-statistics. The distribution of a decoupled $U$-statistic is not changing if the sequences of random variables put in some coordinates of its kernel function are replaced by an independent copy, and this is a very useful property in the application of symmetrization arguments. Beside this, the usual arguments applied in calculation with usual $U$-statistics can be adapted to the study of decoupled $U$-statistics. Now I formulate the following version of Proposition 7.1.

**Proposition 8.2.** *Consider a class of functions $f \in \mathcal{F}$ on the $k$-fold product $(X^k, \mathcal{X}^k)$ of a measurable space $(X, \mathcal{X})$, a probability measure $\mu$ on $(X, \mathcal{X})$ together with a sequence of independent and $\mu$ distributed random variables $\xi_1, \ldots, \xi_n$ which satisfy the conditions of Proposition 7.1. Let us take $k$ independent copies $\xi_1^{(j)}, \ldots, \xi_n^{(j)}$, $1 \leq j \leq k$, of the random sequence $\xi_1, \ldots, \xi_n$, and consider the decoupled $U$-statistics $\bar{I}_{n,k}(f)$, $f \in \mathcal{F}$, defined with their help by formula (53). There exists a sufficiently large constant $K = K(k)$ together with some number $\gamma = \gamma(k) > 0$ and threshold index $A_0 = A_0(k) > 0$ depending only on the order $k$ of the decoupled $U$-statistics we consider such that if $n\sigma^2 > K(L+\beta) \log n$ with $\beta = \max\left(\frac{\log D}{\log n}, 0\right)$, then the (degenerate) decoupled $U$-statistics $\bar{I}_{n,k}(f)$, $f \in \mathcal{F}$, satisfy the following version of inequality (47):*

$$P\left(\sup_{f \in \mathcal{F}} |n^{-k/2}\bar{I}_{n,k}(f)| \geq An^{k/2}\sigma^{k+1}\right) \leq e^{-\gamma A^{1/2k}n\sigma^2} \quad \text{if } A \geq A_0. \tag{57}$$

Proposition 8.2 and Theorem 8.1 imply Proposition 7.1. Hence it is enough to concentrate on the proof of Proposition 8.2. It is natural to try to adapt the method applied in the proof of Proposition 7.1 in the case $k = 1$. I try to explain what kind of new problems appear in the multivariate case and how to overcome them.

The proof of Proposition 7.1 was based on a symmetrization type result formulated in Lemma 7.2 and Hoeffding's inequality Theorem 7.3. We have to find the multivariate versions of these results. It is not difficult to find the multivariate version of Hoeffding's inequality. Such a result can be found in [22] Theorem 12.3, or [21] contains an improved version with optimal constant in the exponent. Here I do not formulate this result, I only explain its main content. Let us consider a homogeneous polynomial of Rademacher functions of order $k$. The multivariate version of Hoeffding's inequality states that its tail distribution can be bounded by that of $K\sigma\eta^k$ with some constant $K = K(k)$ depending only on the order $k$ of the homogeneous polynomial, where $\eta$ is a standard normal random variable, and $\sigma^2$ is the variance of the random homogeneous polynomial.

The problem about the multivariate generalization of Lemma 7.2 is much harder. We want to prove the following multivariate version of this result.

**Lemma 8.3.** *Let $\mathcal{F}$ be a class of functions on the space $(X^k, \mathcal{X}^k)$ which satisfies the conditions of Proposition 7.1 with some probability measure $\mu$. Let us have $k$ independent copies $\xi_1^{(j)}, \ldots, \xi_n^{(j)}$, $1 \leq j \leq k$, of a sequence of independent $\mu$ distributed random variables $\xi_1, \ldots, \xi_n$ and a sequence of independent random variables $\varepsilon = (\varepsilon_1, \ldots, \varepsilon_n)$, $P(\varepsilon_l = 1) = P(\varepsilon_l = -1) = \frac{1}{2}$, $1 \leq l \leq n$, which is independent also of the random sequences $\xi_1^{(j)}, \ldots, \xi_n^{(j)}$, $1 \leq j \leq k$. Consider the decoupled U-statistics $\bar{I}_{n,k}(f)$ defined with the help of these random variables by formula (53) together with their randomized version $\bar{I}_{n,k}^\varepsilon(f)$ defined in (54) for all $f \in \mathcal{F}$. There exists some constant $A_0 = A_0(k) > 0$ such that the inequality*

$$P\left(\sup_{f \in \mathcal{F}} n^{-k/2}\left|\bar{I}_{n,k}(f)\right| > An^{k/2}\sigma^{k+1}\right) \tag{58}$$

$$< 2^{k+1}P\left(\sup_{f \in \mathcal{F}}\left|\bar{I}_{n,k}^\varepsilon(f)\right| > 2^{-(k+1)}An^k\sigma^{k+1}\right) + Bn^{k-1}e^{-A^{1/(2k-1)}n\sigma^2/k}$$

*holds for all $A \geq A_0$ with some appropriate constant $B = B(k)$. One can choose for instance $B = 2^k$ in this result.*

The estimate (58) in Lemma 8.3 is similar to formula (48) in Lemma 7.2. There is a slight difference between them, because the right-hand side of (58) contains an additional constant term. But this term is sufficiently small, and its presence causes no problem as we try to prove Proposition 8.2 by means of Lemma 8.3. In this proof we want to estimate the distribution of the supremum of the decoupled U-statistics $\bar{I}_{n,k}(f)$, $f \in \mathcal{F}$, defined in formula (53), and Lemma 8.3 helps us in reducing this problem to an analogous one, where these decoupled U-statistics are replaced by the randomized decoupled U-statistics $\bar{I}_{n,k}^\varepsilon(f)$, defined in formula (54). This reduced problem can be studied by taking the conditional probability of the event whose probability is considered at the right-hand side of (58) with respect to the condition that all random variables $\xi_l^{(j)}$, $1 \leq j \leq k$, $1 \leq l \leq n$, take a prescribed value. These conditional probabilities can be estimated by means of the multivariate version of the Hoeffding inequality, and then an adaptation of the method described in the previous section supplies the proof of Proposition 8.2. The proof is harder in this new case, but no new principal difficulty arises.

Lemma 7.2 was proved by means of a simple result formulated in Lemma 7.4 which enabled us to introduce the randomizing terms $\varepsilon_j$, $1 \leq j \leq n$. In this result we have taken beside the original sequence $\xi_1, \ldots, \xi_n$ an independent copy $\bar\xi_1, \ldots, \bar\xi_n$. In the next Lemma 8.4 I formulate a multivariate version of Lemma 7.4 which may help in the proof of Lemma 8.3. In its formulation I introduce beside the $k$ independent copies $\xi_1^{(j)}, \ldots, \xi_n^{(j)}$, $1 \leq j \leq k$, of the original sequence of independent, identically distributed random variables $\xi_1, \ldots, \xi_n$ appearing in the definition of a decoupled $U$-statistic of order $k$ another $k$ independent copies $\bar\xi_1^{(j)}, \ldots, \bar\xi_n^{(j)}$, $1 \leq j \leq k$, of this sequence. Because of notational convenience I reindex them, and I shall deal in Lemma 8.4 with $2k$ independent copies $\xi_1^{(j,1)}, \ldots, \xi_n^{(j,1)}$ and $\xi_1^{(j,-1)}, \ldots, \xi_n^{(j,-1)}$, $1 \leq j \leq k$, of the original sequence $\xi_1, \ldots, \xi_n$.

Now I formulate Lemma 8.4.

**Lemma 8.4.** *Let us have a (non-empty) class of functions $\mathcal{F}$ of $k$ variables $f(x_1, \ldots, x_k)$ on a measurable space $(X^k, \mathcal{X}^k)$ together with $2k$ independent copies $\xi_1^{(j,1)}, \ldots, \xi_n^{(j,1)}$ and $\xi_1^{(j,-1)}, \ldots, \xi_n^{(j,-1)}$, $1 \leq j \leq k$, of a sequence of independent and identically distributed random variables $\xi_1, \ldots, \xi_n$ on $(X, \mathcal{X})$ and another sequence of independent random variables $\varepsilon_1, \ldots, \varepsilon_n$, $P(\varepsilon_j = 1) = P(\varepsilon_j = -1) = \frac12$, $1 \leq j \leq n$, independent of all previously considered random sequences. Let us denote the class of sequences of length $k$ consisting of $\pm 1$ digits by $V_k$, and let $m(v)$ denote the number of digits $-1$ in a sequence $v = (v(1), \ldots, v(k)) \in V_k$. Let us introduce with the help of the above notations the random variables $\tilde I_{n,k}(f)$ and $\tilde I_{n,k}(f, \varepsilon)$ as*

$$\tilde I_{n,k}(f) = \frac{1}{k!} \sum_{v \in V_k} (-1)^{m(v)} \sum_{\substack{1 \leq l_r \leq n, \ r=1,\ldots,k \\ l_r \neq l_{r'} \ if \ r \neq r'}} f\left(\xi_{l_1}^{(1,v(1))}, \ldots, \xi_{l_k}^{(k,v(k))}\right) \qquad (59)$$

*and*

$$\tilde I_{n,k}(f, \varepsilon) = \frac{1}{k!} \sum_{v \in V_k} (-1)^{m(v)} \sum_{\substack{1 \leq l_r \leq n, \ r=1,\ldots,k \\ l_r \neq l_{r'} \ if \ r \neq r'}} \varepsilon_{l_1} \cdots \varepsilon_{l_k} f\left(\xi_{l_1}^{(1,v(1))}, \ldots, \xi_{l_k}^{(k,v(k))}\right)$$
$$(60)$$

*for all $f \in \mathcal{F}$. The joint distributions of the random variables $\{\tilde I_{n,k}(f); \ f \in \mathcal{F}\}$ and $\{\tilde I_{n,k}(f, \ \varepsilon); f \in \mathcal{F}\}$ defined in formulas (59) and (60) agree.*

The proof of Lemma 8.4 can be found as Lemma 11.5 in [22]. Actually, this proof is not difficult. Let us observe that the inner sum in formula (59) is a decoupled $U$-statistic, and in formula (60) it is a randomized decoupled $U$-statistic. (Actually they are multiplied by $k!$). In formulas (59) and (60) such a linear combination of these expressions was taken which is similar to the formula appearing in the definition of Stieltjes measures.

Let us list the functions in the class of functions $\mathcal{F}$ in Lemma 8.3 in the form $\{f_1, f_2, \dots\} = \mathcal{F}$, and introduce the quantities

$$Z_p = \frac{n^{-k/2}}{k!} \sum_{\substack{1 \le l_r \le n, \ r=1,\dots,k \\ l_r \ne l_{r'}, \text{ if } r \ne r'}} f_p\left(\xi_{l_1}^{(1,1)}, \dots, \xi_{l_k}^{(k,1)}\right), \quad p = 1, 2, \dots, \quad (61)$$

and

$$\bar{Z}_p = Z_p - n^{-k/2} \tilde{I}_{n,k}(f_p), \quad p = 1, 2, \dots, \quad (62)$$

with the random variables $\tilde{I}_{n,k}(f)$ introduced in (59) with the function $f = f_p$. We would like to prove Lemma 8.3 with the help of Lemma 8.4. This can be done with the help of some calculations, but this requires to overcome some very hard problems. We should like to bound a probability of the form $P\left(\sup_{1 \le p < \infty} Z_p > u\right)$ from above with the help of a probability of the form $P\left(\sup_{1 \le p < \infty} (Z_p - \bar{Z}_p) > \frac{u}{2}\right)$ for all sufficiently large numbers $u$. The question arises how to prove such an estimate. This problem is the most difficult part of the proof.

In the case $k = 1$ considered in the previous section the analogous problem could be simply solved by means of a Symmetrization Lemma formulated in Lemma 7.5. This Lemma cannot be applied in the present case, because it has an important condition, it demands that the sequences of random variables $Z_p$, $p = 1, 2, \dots$, and $\bar{Z}_p$, $p = 1, 2, \dots$, should be independent. In the problem of Section 7 we could work with such sequences which satisfy this condition. On the other hand, the sequences $Z_p$ and $\bar{Z}_p$, $p = 1, 2, \dots$, defined in formulas (61) and (62) we have to work with now are not independent in the case $k \ge 2$. They satisfy some weak sort of independence, and the problem is how to exploit this to get the estimates we need.

Let us first formulate such a version of the Symmetrization Lemma which can be applied also in the problem investigated now. This is done in the next Lemma 8.5.

**Lemma 8.5 (Generalized version of the Symmetrization Lemma).** *Let $Z_p$ and $\bar{Z}_p$, $p = 1, 2, \dots$, be two sequences of random variables on a probability space $(\Omega, \mathcal{A}, P)$. Let a $\sigma$-algebra $\mathcal{B} \subset \mathcal{A}$ be given on the probability space $(\Omega, \mathcal{A}, P)$ together with a $\mathcal{B}$-measurable set $B$ and two numbers $\alpha > 0$ and $\beta > 0$ such that the random variables $Z_p$, $p = 1, 2, \dots$, are $\mathcal{B}$ measurable, and the inequality*

$$P(|\bar{Z}_p| \le \alpha | \mathcal{B})(\omega) \ge \beta \quad \text{for all } p = 1, 2, \dots \text{ if } \omega \in B \quad (63)$$

*holds. Then*

$$P\left(\sup_{1 \le p < \infty} |Z_p| > \alpha + u\right) \le \frac{1}{\beta} P\left(\sup_{1 \le p < \infty} |Z_p - \bar{Z}_p| > u\right) + (1 - P(B)) \quad (64)$$

*for all $u > 0$.*

The proof of Lemma 8.5 is contained together with its proof in [22] under the name Lemma 13.1, and the proof is not hard. It consists of a natural adaptation of the proof of the original Symmetrization Lemma, presented in Lemma 7.5. The hard problem is to check the condition in formula (63) in concrete applications. In our case we would like to apply this lemma to the random variables $Z_p$ and $\bar{Z}_p$, $p = 1, 2, \ldots$, defined in formulas (61) and (62) together with the $\sigma$-algebra $\mathcal{B} = \mathcal{B}(\xi_1^{(j,1)}, \ldots, \xi_n^{(j,1)}, 1 \leq j \leq k)$ generated by the random variables $\xi_1^{(j,1)}, \ldots, \xi_n^{(j,1)}$, $1 \leq j \leq k$. We would like to show that relation (63) holds with this choice on a set $B$ of probability almost 1. (Let me emphasize that in (63) a set of inequalities must hold for all $p = 1, 2, \ldots$ simultaneously if $\omega \in B$.)

In the analogous problem considered in Section 7 condition (51) had to be checked with some appropriate constants $\alpha > 0$ and $\beta > 0$ for the random variables $\bar{Z}_p$, $p = 1, 2, \ldots$, defined in formula (52). This could be done fairly simply by the calculation of the variance of the random variables $\bar{Z}_p$, $p = 1, 2, \ldots$. A natural adaptation of this approach is to bound from above the supremum $\sup_{1 \leq p < \infty} E\left(\bar{Z}_p^2 | \mathcal{B}\right)$ of the conditional second moments of the random variables $\bar{Z}_p$, $1 \leq p < \infty$, defined in (62) with respect to the $\sigma$-algebra $\mathcal{B}$ and to show that this expression is small with large probability. I have followed this approach in [19] and [22]. One can get the desired estimates, but many unpleasant technical details have to be tackled in the proof. I do not discuss here all details, I only briefly explain what kind of problems we meet when try to apply this method in the special case $k = 2$ and give some indications how they can be overcome.

In the case $k = 2$ the definition of $\bar{Z}_p$ is very similar to that of $n^{-k/2}\bar{I}_{n,2}(f_p)$ defined in (59) with the function $f = f_p$. The only difference is that in the definition of $Z_p$ we have to take the values $v = (1, -1)$, $v = (-1, 1)$ and $v = (-1, -1)$ in the outer sum, i.e. the term $v = (1, 1)$ is dropped, and we multiply by $(-1)^{m(v)+1}$ instead of $(-1)^{m(v)}$. We can get the desired estimate on the conditional supremum of second moments if we can prove a good estimate on the conditional second moments of the supremum of the inner sums in $\bar{I}_{n,2}(f_p)$, $1 \leq p < \infty$, in the case of each index $v = (1, -1)$, $v = (-1, 1)$ and $v = (-1, -1)$. If we can get a good estimate in the case $v = (1, -1)$, then we can get it in the remaining cases, too. So we have to give a good bound on the expression

$$\sup_{1 \leq p < \infty} E\left(\frac{1}{n}\left(\sum_{1 \leq l_r \leq n, \; r=1,2, \; l_1 \neq l_2} f_p\left(\xi_{l_1}^{(1,1)}, \xi_{l_2}^{(2,-1)}\right)\right)^2 \middle| \mathcal{B}\right). \qquad (65)$$

Moreover, since the sequence of random variables $\xi_l^{(2,-1)}$, $1 \leq l \leq n$, is independent of the $\sigma$-algebra $\mathcal{B}$, and the canonical property of the functions $f_p$ implies some orthogonalities, the estimation of the expression in (65) can be simplified. A detailed calculation shows that it is enough to prove the following inequality:

Let us have a countable class $\mathcal{F}$ of canonical functions $f(x, y)$ with respect to a probability measure $\mu$ on the second power $(X^2, \mathcal{X}^2)$ of a measurable space $(X, \mathcal{X})$, which is $L_2$-dense with some exponent $L$ and parameter $D$, (the probability measure $\mu$ is living in the space $(X, \mathcal{X})$) together with a sequence of independent and $\mu$-distributed random variables $\xi_1, \ldots, \xi_n$, $n \geq 2$, on $(X, \mathcal{X})$, and let the relations

$$\int f(x, y)^2 \mu(dx)\mu(dy) \leq \sigma^2, \quad \sup |f(x, y)| \leq 1 \qquad \text{for all } f \in \mathcal{F}$$

hold with some number $0 < \sigma^2 \leq 1$ which satisfies the relation $n\sigma^2 \geq K(L + \beta)\log n$ with $\beta = \max\left(\frac{\log D}{\log n}, 0\right)$ and a sufficiently large fixed constant $K > 0$. Then the inequality

$$P\left(\sup_{f \in \mathcal{F}} \frac{1}{n} \int \left(\sum_{l=1}^{n} f(\xi_l, y)\right)^2 \mu(dy) \geq A^2 n\sigma^4\right) \leq \exp\left\{-A^{1/3} n\sigma^2\right\} \qquad (66)$$

holds if $A \geq A_0$ with some sufficiently large fixed constant $A_0$.

Inequality (66) is similar to relation (47) in Proposition 7.1 in the case $k = 1$, but it does not follow from it. (It follows from (47) in the special case when the function $f$ does not depend on the argument $y$ with respect to which we integrate.) On the other hand, inequality (66) can be proved by working out a similar, although somewhat more complicated symmetrization argument and induction procedure as it was done in the proof of Proposition 7.1 in the case $k = 1$. After this, inequality (66) enables us to work out the symmetrization argument we need to prove Proposition 7.1 for $k = 2$. This procedure can be continued for all $k = 2, 3, \ldots$. If we have already proved Proposition 7.1 for some $k$, then an inequality can be formulated and proved with the help of the already known results which enable us to carry out that symmetrization procedure which is needed in the proof of Proposition 7.1 in the case $k + 1$. This is a rather cumbersome method with a lot of technical details, hence its detailed explanation had to be omitted from an overview paper. In the work [22] Sections 13, 14 and 15 deal only with the proof of Proposition 7.1. Section 13 contains the proof of some preparatory results and the formulation of the inductive statements we have to prove to get the result of Proposition 7.1, Section 14 contains the proof of the Symmetrization arguments we need, and finally the proof is completed with their help in Section 15.

There is an interesting theory of Talagrand about so-called concentration inequalities. This theory has some relation to the questions discussed in this paper. In the last section this relation will be discussed together with some other results and open problems.

## 9. Relation with other results and some open problems

Talagrand worked out a deep theory about so-called concentration inequalities. (See his overview in paper [33] about this subject.) His results are closely related to the supremum estimates described in this paper. First I discuss this relation.

### 9.1. On Talagrand's concentration inequalities

Talagrand considered a sequence of independent random variables $\xi_1, \ldots, \xi_n$, a class of functions $\mathcal{F}$, took the partial sums $\sum_{j=1}^{n} f(\xi_j)$ for all functions $f \in \mathcal{F}$, and investigated their supremum. He proved such estimates which state that this supremum is very close to its expected value, (it is concentrated around it). The following theorem in paper [34] is a typical result in this direction.

**Theorem 9.1 (Theorem of Talagrand).** *Consider $n$ independent and identically distributed random variables $\xi_1, \ldots, \xi_n$ with values in some measurable space $(X, \mathcal{X})$. Let $\mathcal{F}$ be some countable family of real-valued measurable functions of $(X, \mathcal{X})$ such that $\|f\|_\infty \leq b < \infty$ for every $f \in \mathcal{F}$. Let $Z = \sup_{f \in \mathcal{F}} \sum_{i=1}^{n} f(\xi_i)$ and $v = E(\sup_{f \in \mathcal{F}} \sum_{i=1}^{n} f^2(\xi_i))$. Then for every positive number $x$,*

$$P(Z \geq EZ + x) \leq K \exp\left\{-\frac{1}{K'}\frac{x}{b}\log\left(1 + \frac{xb}{v}\right)\right\} \tag{67}$$

*and*

$$P(Z \geq EZ + x) \leq K \exp\left\{-\frac{x^2}{2(c_1 v + c_2 bx)}\right\}, \tag{68}$$

*where $K$, $K'$, $c_1$ and $c_2$ are universal positive constants. Moreover, the same inequalities hold when replacing $Z$ by $-Z$.*

Inequality (67) can be considered as a generalization of Bennett's inequality, inequality (68) as a generalization of Bernstein's inequality. In these estimates the distribution of the supremum of possibly infinitely many partial sums of independent and identically distributed functions are considered. A remarkable feature of Theorem 9.1 is that it imposes no condition about the structure of the class of functions $\mathcal{F}$. In this respect it differs from Theorems 6.2 and 6.3 in this paper, where such a class of functions $\mathcal{F}$ is considered which satisfies a so-called $L_2$-density property.

Talagrand's study was also continued by other authors who got interesting results. In particular, the works of M. Ledoux [16] and P. Massart [26] are worth mentioning. In these works the above mentioned result was improved. Such a version was proved which also holds for the supremum of appropriate classes of sums of independent but not necessarily identically distributed random variables. (On the other hand, I do not know of such a generalization in which $U$-statistics of higher order are considered.) The improvements of these works consist for instance in a version of Theorem 9.1 where the quantity $v = E(\sup_{f \in \mathcal{F}} \sum_{i=1}^{n} f^2(\xi_i))$ is replaced by $\sigma^2 = \sup_{f \in \mathcal{F}} \sum_{i=1}^{n} \text{Var}(f(\xi_i))$, i.e. the supremum of the expectation of the individual partial sums $\sum_{i=1}^{n} f^2(\xi_i)$ is considered (the statement that $\sigma^2$ equals the supremum of the expected values of the partial sums $\sum_{i=1}^{n} f^2(\xi_i)$ holds if $Ef(\xi_i) = 0$ for all random variables $\xi_i$ and functions $f$) instead of the second moment of the supremum of these partial sums.

On the other hand, the estimates in Theorem 9.1 contain the expected value $EZ = E\left(\sup_{f \in \mathcal{F}} \sum_{i=1}^{n} f(\xi_i)\right)$, and this quantity appears in all concentration type inequalities. This fact has deep consequences which deserve a more detailed discussion.

Let us consider Theorem 9.1 or one of its improvements and try to understand what kind of solution they provide for problem b) or b') formulated in Section 1 in the case $k = 1$. They supply a good estimate on the probabilities we consider for the numbers $u \geq n^{-1/2} EZ = n^{-1/2} E(\sup_{f \in \mathcal{F}} \sum_{i=1}^{n} f(\xi_i))$. But to apply these results we need a good estimate on the expectation $EZ$ of the supremum of the partial sums we consider, and the proof of such an estimate is a highly non-trivial problem.

Let us consider problem b') (in the case $k = 1$) for such a class of functions $\mathcal{F}$ which satisfies the conditions of Theorem 6.3. The considerations taken in Section 6 show that there are such classes of functions $\mathcal{F}$ which satisfy the conditions of Theorem 6.3, and for which the probability $P(\sup_{f \in \mathcal{F}} n^{-1/2} \sum_{i=1}^{n} f(\xi_i) > \alpha \sigma \log \frac{2}{\sigma})$ is almost 1 with an appropriate small number $\alpha > 0$ for all large enough sample sizes $n$. (Here the number $\sigma$ is the same as in Theorem 6.3.) This means that $En^{-1/2}Z \geq (\alpha - \varepsilon)\sigma \log \frac{2}{\sigma}$ for all $\varepsilon > 0$ if the sample size $n$ of the sequence $\xi_1, \ldots, \xi_n$ is greater than $n_0 = n_0(\varepsilon, \sigma)$. Some calculation also shows that under the conditions of Theorem 6.3 $En^{-1/2}Z \leq K\sigma \log \frac{2}{\sigma}$ with an appropriate number $K > 0$. (In this calculation some difficulty may arise, because Theorem 6.3 for $k = 1$ does not yield a good estimate if $u \geq \sqrt{n}\sigma^2$. But we can write $P(\sup_{f \in \mathcal{F}} n^{-1/2} \sum_{i=1}^{n} f(\xi_i) > u) \leq e^{-\alpha(u/\bar{\sigma})^2} = e^{-\alpha u \sqrt{n}}$ with $\bar{\sigma}^2 = un^{-1/2}$ if $u \geq \sqrt{n}\sigma^2$, and this estimate is sufficient for us. We get the upper bound we formulated for $n^{-1/2}EZ$ from Theorem 6.3 only under the condition $n\sigma^2 \geq \text{const.} \log \frac{2}{\sigma}$ with some appropriate constant. It can be seen that this condition is really needed, it appeared not because of the weakness of our method. I omit the details of the calculation.) Then the concentration inequality Theorem 9.1, or more precisely its improvement, Theorem 3 in paper [26] which gives a similar inequality, but with the quantity $\sigma^2$ instead of $v$ implies Theorem 6.3 in the case $k = 1$. This means that Theorem 6.3 can be deduced from concentration type inequalities in the case $k = 1$ if we can show that under its conditions $En^{-1/2}Z \leq K\sigma \log \frac{2}{\sigma}$ with some appropriate $K > 0$ depending only on the exponent and parameter of the $L_2$-dense class $\mathcal{F}$. Such an estimate can be proved (see the proof in [8] on the basis of paper [33]), but it requires rather long and non-trivial considerations. I prefer a direct proof of Theorem 6.3.

Finally I discuss a refinement of Theorems 4.1 and 4.3 promised in a remark at the end of Section 4 together with some open problems.

### 9.2. Some refinements of the estimate in Theorems 4.1 and 4.3

If we have a bound on the $L_2$ and $L_\infty$ norm of the kernel function $f$ of a U-statistic $I_{k,n}(f)$, but we have no additional information about the behaviour of $f$, (and such a situation is quite common in mathematical statistics problems), then the estimate of Theorem 4.3 about the distribution of U-statistics cannot be considerably improved. On the other hand, one would like to prove such a multi-dimensional type version of the large deviation theorem about partial sums of independent random variables which gives a good asymptotic formula for the probability $P(n^{-k/2}I_{k,n}(f) > u)$ for large values $u$. Such an estimate should depend on the function $f$. A similar question can be posed about the distribution of multiple Wiener-Itô integrals $Z_{n,k}(f)$ if $k \geq 2$, because the distribution of such random integrals (unlike the degenerate case $k = 1$) is not determined by their variance.

Such large deviation problems are very hard, and I know of no result in this direction. On the other hand, some quantities can be introduced which enable us to give a better estimate on the distribution of Wiener–Itô integrals or U-statistics in the case of their knowledge. Such results were known for Wiener–Itô integrals $Z_{\mu,2}(f)$ and U-statistics $I_{n,2}(f)$ of order 2 earlier, and quite recently they were generalized for all $k \geq 2$. I describe them and show that they are useful in the solution of some problems. My formulation will differ a little bit from the previous ones. In particular, I shall speak about Wiener–Itô integrals where previous authors considered only polynomials of Gaussian random vectors. But the Wiener–Itô integral presentation of these results seems to be more natural for me. First I formulate the estimate about Wiener–Itô integrals of order 2 proved in [12] by Hanson and Wright.

**Theorem 9.2.** *Let a two-fold Wiener–Itô integral*

$$Z_{\mu,2}(f) = \int f(x,y)\mu_W(dx)\mu_W(dy)$$

*be given, where $\mu_W$ is a white noise with a non-atomic reference measure $\mu$, and the function $f$ satisfies the inequalities*

$$\int f(x,y)^2 \mu(dx)\mu(dy) \leq \sigma^2 \tag{69}$$

*and*

$$\int f(x,y)g_1(x)g_2(y)\mu(dx)\mu(dy) \leq D \tag{70}$$

*with some number $D > 0$ for all functions $g_1$ and $g_2$ such that $\int g_j^2(x)\mu(dx) \leq 1$, $j = 1, 2$. There exists a universal constant $K > 0$ such that the inequality*

$$P(|Z_{\mu,2}(f)| > u) \leq K \exp\left\{-\frac{1}{K}\min\left(\frac{u^2}{\sigma^2}, \frac{u}{D}\right)\right\} \tag{71}$$

*holds for all $u > 0$.*

As it was remarked in Section 4 we can assume without violating the generality that the function $f$ in the definition of Wiener–Itô integrals is symmetric. In this case Theorem 9.2 can be reformulated to a simpler statement.

To do this let us define with the help of the (symmetric) function $f$ the following so-called Hilbert–Schmidt operator $A_f$ in the $L_2(\mu)$ space of square integrable functions with respect to the measure $\mu$: $A_f v(x) = \int f(x,y)v(y)\mu(dy)$ for all $L_2(\mu)$ measurable functions $v(\cdot)$. It is known that $A_f$ is a compact, self-adjoint operator, hence it has a discrete spectrum. Let $\lambda_1, \lambda_2, \ldots$ denote the eigenvalues of the operator $A_f$. It follows from the theory of Hilbert-Schmidt operators and the Itô formula for multiple Wiener–Itô integrals that the identity $Z_{\mu,2}(f) = \sum_{j=1}^{\infty} \lambda_j(\eta_j^2 - 1)$ holds with some appropriately defined independent standard normal random variables $\eta_1, \eta_2, \ldots$. Beside this, $\sum_{j=1}^{\infty} \lambda_j^2 = \int f^2(x,y)\mu(dx)\mu(dy)$. Hence condition (69) can be reformulated as $\sum_{j=1}^{\infty} \lambda_j^2 \leq \sigma^2$, and condition (70) is equivalent to the statement that $\sup_j |\lambda_j| \leq D$. In such a way Theorem 9.2 can be reduced to another statement whose proof is simpler.

Theorem 9.2 yields a useful estimate if $D^2 \ll \sigma^2$. In this case it states that for large numbers $u$ the bound $P(Z_{\mu,2}(f) > u) \leq \text{const.} e^{-u/2\sigma}$ supplied by Theorem 4.1 can be improved to the bound $P(Z_{\mu,2}(f) > u) \leq \text{const.} e^{-u/KD}$. The correction term $\frac{u^2}{\sigma^2}$ at the right-hand side of (71) is needed to get an estimate which holds for all $u > 0$. It may be worthwhile recalling the following result (see [27] or [17], Theorem 6.6). All $k$-fold Wiener–Itô integrals $Z_{\mu,k}(f)$ satisfy the inequality $P(|Z_{\mu,k}(f)| > u) > K e^{-Au^{2/k}}$ with some $K = K(f,\mu) > 0$ and $A = A(f,\mu) > 0$. There is a strictly positive number $A = A(f,\mu)$ in the exponent of the last relation, but the proof of [27] yields no explicit lower bound for it.

There is a similar estimate about the distribution of degenerate $U$-statistics of order 2. This is the content of the following Theorem 9.3.

**Theorem 9.3.** *Let a sequence $\xi_1, \ldots, \xi_n$ of independent $\mu$ distributed random variables be given together with a function $f(x,y)$ canonical with respect to the measure $\mu$, and consider the (degenerate) $U$-statistic $I_{n,2}(f)$ defined in (3) with the help of the above quantities. Let us assume that the function $f$ satisfies conditions (69) and (70) with some $\sigma > 0$ and $D > 0$, and also the relations*

$$\sup_x \int f^2(x,y)\mu(dy) \leq A_1, \quad \sup_y \int f^2(x,y)\mu(dx) \leq A_2, \quad \sup_{x,y} |f(x,y)| \leq B \tag{72}$$

*hold with some appropriate constants $A_1 > 0$, $A_2 > 0$ and $B > 0$. Then there exists a universal constant $K > 0$ such that the inequality*

$$P\left(n^{-1}|I_{n,2}| > u\right) \leq K \exp\left\{-\frac{1}{K}\left(\frac{u^2}{\sigma^2}, \frac{u}{D}, \frac{n^{1/3}u^{2/3}}{(A_1 + A_2)^{1/3}}, \frac{n^{1/2}u^{1/2}}{B^{1/2}}\right)\right\} \tag{73}$$

*is valid for all $u > 0$.*

Theorem 9.3 was proved in [9]. The estimate of Theorem 9.3 is similar to that of Theorem 9.2, the difference between them is that in formula (73) some additional correction terms had to be inserted to make it valid for all $u > 0$. But the proof of Theorem 9.3 is much harder. It can be shown that the estimate (73) implies that of Theorem 4.3 in the special case $k = 2$ if we disregard the appearance of the not explicitly defined universal constant $K$ in it.

To see this observe that Theorem 4.3 contains the conditions $u \leq n^{1/2}\sigma^2$ and $B \leq 1$ which imply that $\frac{n^{1/3}u^{2/3}}{(A_1+A_2)^{1/2}} \geq \frac{1}{\sqrt{2}}\left(\frac{u}{\sigma^2}\right)^{2/3} u^{2/3} = \frac{1}{\sqrt{2}}\left(\frac{u}{\sigma}\right)^{4/3}$, and $\frac{n^{1/2}u^{1/2}}{B^{1/2}} \geq \frac{u}{\sigma^2}u^{1/2} = \sigma^{-1/2}\left(\frac{u}{\sigma}\right)^{3/2} \geq \left(\frac{u}{\sigma}\right)^{3/2}$, since $\sigma \leq 1$ in this case. Beside this, $\frac{u}{D} \geq \frac{u}{\sigma}$. The above relations imply that in the case $u \geq \sigma$ the estimate (73) is weakened if the expression in its exponent is replaced by $\frac{1}{\sqrt{2}K}\frac{u}{\sigma}$. Theorem 4.3 trivially holds if $0 \leq u \leq \sigma$.

Theorem 9.3 is useful in such problems where a refinement of the estimate in Theorem 4.3 is needed which exploits better the properties of the kernel function $f$ of a degenerate $U$-statistics of order 2. Such a situation appears in paper [10], where the law of iterated logarithm is investigated for degenerate $U$-statistics of order 2.

Let us consider an infinite sequence $\xi_1, \xi_2, \ldots$ of independent $\mu$ distributed random variables together with a function $f$ canonical with respect to the measure $\mu$, and define the degenerate $U$-statistic $I_{n,2}(f)$ with their help for all $n = 1, 2, \ldots$. In paper [10] the necessary and sufficient condition of the iterated logarithm is given for such a sequence. More explicitly, it is proved that

$$\limsup_{n\to\infty} \frac{|I_{n,2}(f)|}{n\log\log n} < \infty \quad \text{with probability } 1$$

if and only if the following two conditions are satisfied:

a) $\int_{\{(x,y):\, |f(x,y)|\leq u\}} f^2(x,y)\mu(dx)\mu(dy) \leq C\log\log u$ with some $C < \infty$ for all $u \geq 10$.

b) $\int f(x,y)g(x)h(y)\mu(dx)\mu(dy) \leq C$ with some appropriate $C < \infty$ for all such pairs of functions $g$ and $h$ which satisfy the relations $\int g^2(x)\mu(dx) \leq 1$, $\int h^2(x)\mu(dx) \leq 1$, $\sup_x |g(x)| < \infty$, $\sup_x |h(x)| < \infty$.

The above result is proved by means of a clever truncation of the terms in the $U$-statistics and an application of the estimation of Theorem 9.3 for these truncated $U$-statistics. It has the form one would expect by analogy with the classical law of iterated logarithm for sums of independent, identically distributed random variables with expectation zero, but it also has an interesting, unexpected feature. The classical law of iterated logarithm for sums of iid. random variables holds if and only if the terms in the sum have finite variance. (The only if part is proved in paper [7] or [30].) The above formulated law of iterated logarithm for degenerate $U$-statistics also holds in the case of finite second moment, i.e. if $Ef^2(\xi_1,\xi_2)<\infty$, but as the authors in [10] show in an example, there are also cases when it holds, although $Ef^2(\xi_1,\xi_2)=\infty$. Paper [11] is another example where Theorem 9.3 can be successfully applied to solve certain problems.

To formulate the generalization of Theorems 9.2 and 9.3 for general $k\geq2$ some notations have to be introduced. Given a finite set $A$ let $\mathcal{P}(A)$ denote the set of all its partitions. If a partition $P=\{B_1,\ldots,B_s\}\in\mathcal{P}(A)$ consists of $s$ elements then we say that this partition has order $s$, and write $|P|=s$. In the special case $A=\{1,\ldots,k\}$ the notation $\mathcal{P}(A)=\mathcal{P}_k$ will be used. Given a measurable space $(X,\mathcal{X})$ with a probability measure $\mu$ on it together with a finite set $B=\{b_1,\ldots,b_j\}$ let us introduce the following notations. Take $j$ different copies $(X_{b_r},\mathcal{X}_{b_r})$ and $\mu_{b_r}$, $1\leq r\leq j$, of this measurable space and probability measure indexed by the elements of the set $B$, and define their product $(X^{(B)},\mathcal{X}^{(B)},\mu^{(B)})=\left(\prod\limits_{r=1}^{j}X_{b_r},\prod\limits_{r=1}^{j}\mathcal{X}_{b_r},\prod\limits_{r=1}^{j}\mu_{b_r}\right)$. The points $(x_{b_1},\ldots,x_{b_j})\in X^{(B)}$ will be denoted by $x^{(B)}\in X^{(B)}$ in the sequel. With the help of the above notations I introduce the quantities needed in the formulation of the generalization of Theorems 9.2 and 9.3.

Let a function $f=f(x_1,\ldots,x_k)$ be given on the $k$-fold product $(X^k,\mathcal{X}^k,\mu^k)$ of a measurable space $(X,\mathcal{X})$ with a probability measure $\mu$. For all partitions $P=\{B_1,\ldots,B_s\}\in\mathcal{P}_k$ of the set $\{1,\ldots,k\}$ consider the functions $g_r\left(x^{(B_r)}\right)$ on the space $X^{(B_r)}$, $1\leq r\leq s$, and define with their help the quantity

$$\alpha(P)=\alpha(P,f,\mu) \tag{74}$$

$$=\sup_{g_1,\ldots,g_s}\left\{\int f(x_1,\ldots,x_k)g_1\left(x^{(B_1)}\right)\cdots g_s\left(x^{(B_s)}\right)\mu(dx_1)\ldots\mu(dx_k):\right.$$

$$\left.\int g_r^2\left(x^{(B_r)}\right)\mu^{(B_r)}\left(dx^{(B_r)}\right)\leq1\quad\text{for all }1\leq r\leq s\right\}.$$

In the estimation of Wiener–Itô integrals of order $k$ the quantities $\alpha(P)$, $P\in\mathcal{P}$, play such a role as the numbers $D$ and $\sigma^2$ introduced in formulas (69) and (70) in Theorem 9.2. Observe that in the case $|P|=1$, i.e. if $P=\{1,\ldots,k\}$ the identity $\alpha^2(P)=\int f^2(x_1,\ldots,x_k)\mu(dx_1)\ldots\mu(dx_k)$ holds. The following estimate is valid for Wiener–Itô integrals of general order (see [15]).

**Theorem 9.4.** *Let a $k$-fold Wiener–Itô integral $I_{\mu,k}(f)$, $k\geq1$, be defined with the help of a white noise $\mu_W$ with a non-atomic reference measure $\mu$ and a kernel function $f$ of $k$-variable such that $\int f^2(x_1,\ldots,x_k)\mu(dx_1)\ldots\mu(dx_k)<\infty$. There is some universal constant $C(k)<\infty$ depending only of the order $k$ of the random integral such that the inequality*

$$P(|Z_{\mu,k}(f)|>u)\leq C(k)\exp\left\{-\frac{1}{C(k)}\min_{1\leq s\leq k}\min_{P\in\mathcal{P}_k,\,|P|=s}\left(\frac{u}{\alpha(P)}\right)^{2/s}\right\} \tag{75}$$

*holds for all $u>0$ with the quantities $\alpha(P)$, $P\in\mathcal{P}_k$, defined in formula (74).*

Also the following converse estimate holds which shows that the above estimate is sharp. (See again paper [15].) This estimate also yields an improvement of the result in [27] mentioned in this subsection.

**Theorem 9.4′.** *The random integral $Z_{\mu,k}(f)$ considered in Theorem 9.4 also satisfies the inequality*

$$P(|Z_{\mu,k}(f)| > u) \geq \frac{1}{C(k)} \exp\left\{-C(k) \min_{1 \leq s \leq k} \min_{P \in \mathcal{P}_k, |P|=s} \left(\frac{u}{\alpha(P)}\right)^{2/s}\right\}$$

*for all $u > 0$ with some universal constant $C(k) > 0$ depending only on the order $k$ of the integral and the quantities $\alpha(P)$, $P \in \mathcal{P}_k$, defined in formula (74).*

To formulate the result about the distribution of degenerate $U$-statistics for all $k \geq 2$ an analog of the expression $\alpha(P)$ defined in (74) has to be introduced. Let us consider a set $A \subset \{1, \ldots, k\}$ with $|A| = k - r$ elements, $0 \leq r < k$, and a partition $P = \{B_1, \ldots, B_s\} \subset \mathcal{P}(A)$, $1 \leq s \leq k - r$, of $A$ together with a function $f(x_1, \ldots, x_k)$ square integrable with respect to the $k$-fold product $\mu^k$ of a probability measure $\mu$ on a measurable space $(X^k, \mathcal{X}^k)$. Let us introduce, similarly to the definition (74), the quantities

$$\alpha(P, x^{(\{1,\ldots,k\}\setminus A)}) \tag{76}$$
$$= \sup_{g_1,\ldots,g_s}\left\{\int f(x_1,\ldots,x_k)g_1\left(x^{(B_1)}\right)\cdots g_s\left(x^{(B_s)}\right)\mu(dx_1)\ldots\mu(dx_k):\right.$$
$$\left.\int g_u^2\left(x^{(B_u)}\right)\mu^{(B_u)}\left(dx^{(B_u)}\right) \leq 1 \quad \text{for all } 1 \leq u \leq s\right\}$$

for all $x^{(\{1,\ldots,k\}\setminus A)} = \{x_j, j \in \{1,\ldots k,\} \setminus A\} \in X^{\{1,\ldots,k\}\setminus A}$, where $g_u$ are functions defined on $X^{(B_u)}$, $1 \leq u \leq s$, and put

$$\alpha(A, P) = \sup_{x^{\{1,\ldots,k\}\setminus A} \in X^{\{1,\ldots,k\}\setminus A}} \alpha(P, x^{(\{1,\ldots,k\}\setminus A)}). \tag{77}$$

To consider also the case $|A| = k$ when $\{1, \ldots, k\} \setminus A = \emptyset$ let us make the following convention. Let us also speak about the partitions of the empty set by saying that its only partition is the empty set itself. Beside this, put $|\emptyset| = 0$, and

$$\alpha(\{1, \ldots, k\}, \emptyset) = \sup |f(x_1, \ldots, x_k)|. \tag{78}$$

With the help of the above notations the estimate about the distribution of normalized degenerated $U$-statistics proven in [1] can be formulated.

**Theorem 9.5.** *Consider a sequence $\xi_1, \ldots, \xi_n$, $n \geq k$, of independent $\mu$ distributed random variables and a bounded function $f(x_1, \ldots, x_k)$ of $k$, $k \geq 2$, variables canonical with respect to the measure $\mu$. Take the (degenerate) $U$-statistic $I_{n,k}(f)$ defined in (3) with the help of these quantities. There is some* *universal constant $C = C(k) < \infty$ depending only on the order $k$ of this $U$-statistic such that the inequality*

$$P(n^{-k/2}|I_{n,k}(f)| > u)$$
$$\leq C \exp\left\{-\frac{1}{C} \max_{\substack{\{(r,s): \, 0 \leq r < k, \, 1 \leq s \leq k-r\} \\ \cup\{(r,s): \, r=k, \, s=0\}}} \max_{\substack{\{(A,P): \, A \subset \{1,\dots,k\}, \\ |A|=r, \, P \in \mathcal{P}(A), |P|=s\}}} \left(\frac{n^r u^2}{\alpha^2(A,P)}\right)^{1/(2r+s)}\right\}$$

*holds for all $u > 0$ with the above constant $C$ and the quantities $\alpha(A,P)$ defined in (76), (77) and (78).*

It can be seen with the help of some calculation that Theorem 9.5 implies Theorem 4.3 for all orders $k \geq 2$ if we disregard the presence of the unspecified universal constant $C$. (It has to be exploited that under the conditions of Theorem 4.3 $\alpha^2(A,P) \leq \sigma^2$ if $|A| = r$ with $r = 0$, $\alpha(A,P) \leq 1$ for $|A| = r \geq 1$, $\sigma^2 \leq 1$, and $n^{k/2}\sigma^{k+1} \geq u$.)

The proof of Theorems 9.4 and 9.5 is based, similarly to the proof of Theorems 4.1 and 4.3, on a good estimate of the (possibly high) moments of Wiener–Itô integrals and degenerate $U$-statistics. The proofs of these estimates in [1] and [15] are based on many deep and hard inequalities of different authors. One may ask whether the diagram formula, propagated in this work, which gives an explicit formula about these moments cannot be applied in the proof of these results. I think that the answer to this question is in the positive, and even I have some ideas how to carry out such a program. But at the time of writing this work I had not enough time to work out the details.

A natural open problem is to find the large deviation estimates about the tail distribution of multiple Wiener–Itô integrals and $U$-statistics mentioned at the start of this subsection. Such results may better explain why the quantities $\alpha(P)$ and $\alpha(A,P)$ appear in the estimates of Theorems 9.4 and 9.5. It would be interesting to find the true value of the universal constants in these estimates or to get at least some partial results in this direction which would help in solving the following problem:

**Problem.** *Consider a $k$-fold multiple Wiener–Itô integral $Z_{\mu,k}(f)$. Show that its distribution satisfies the relation*

$$\lim_{u \to \infty} u^{2/k} \log P(|Z_{\mu,k}(f)| > u) = K(\mu, f) > 0$$

*with some number $K(\mu, f) > 0$, and determine its value.*

There appear some other natural problems relating to the above results. Thus for instance, it was assumed in all estimates about $U$-statistics discussed in this work that their kernel functions are bounded. A closer study of this condition deserves some attention. It was explained in this paper that its role was to exclude the appearance of some irregular events with relatively large probability which would imply that only weak estimates hold in some cases interesting for us. One may ask whether this condition cannot be replaced by a weaker and more appropriate one in certain problems.

Finally, I mention the following problem.

**Problem.** *Prove an estimate analogous to the result of Theorem 9.5 about the supremum of appropriate classes of U-statistics.*

To solve the above problem one has to tackle some difficulties. In particular, to adapt the method of proof of previous results such a generalization of the multivariate version of Hoeffding's inequality (see [21]) has to be proved about the distribution of homogeneous polynomials of Rademacher functions where the bound depends not only on the variance of these random polynomials, but also on some quantities analogous to the expression $\alpha(P)$ introduced in (74).

## References

[1] ADAMCZAK, R. (2005) Moment inequalities for *U*-statistics. Available at `http://www.arxiv.org/abs/math.PR/0506026`

[2] ALEXANDER, K. (1984) Probability inequalities for empirical processes and a law of the iterated logarithm. *Annals of Probability* **12**, 1041–1067 MR757769

[3] ARCONES, M. A. and GINÉ, E. (1993) Limit theorems for *U*-processes. *Annals of Probability* **21**, 1494–1542 MR1235426

[4] DE LA PEÑA, V. H. and MONTGOMERY–SMITH, S. (1995) Decoupling inequalities for the tail-probabilities of multivariate *U*-statistics. *Annals of Probability*, **23**, 806–816 MR1334173

[5] DUDLEY, R. M. (1998) *Uniform Central Limit Theorems*. Cambridge University Press, Cambridge U.K. MR1720712

[6] DYNKIN, E. B. and MANDELBAUM, A. (1983) Symmetric statistics, Poisson processes and multiple Wiener integrals. *Annals of Statistics* **11**, 739–745 MR707925

[7] FELLER, W. (1968) An extension of the law of the iterated logarithm to variables without variance. *Journal of Mathematics and Mechanics* 343–355 **18** MR233399

[8] GINÉ, E. and GUILLOU, A. (2001) On consistency of kernel density estimators for randomly censored data: Rates holding uniformly over adaptive intervals. *Ann. Inst. Henri Poincaré PR* **37** 503–522 MR1876841

[9] GINÉ, E., LATAŁA, R. and ZINN, J. (2000) Exponential and moment inequalities for *U*-statistics in *High dimensional probability II*. Progress in Probability 47. 13–38. Birkhäuser Boston, Boston, MA. MR1857312

[10] GINÉ, E., KWAPIEŃ, S., LATAŁA, R. and ZINN, J. (2001) The LIL for canonical *U*-statistics of order 2. *Annals of Probability* **29** 520–527 MR1825163

[11] GINÉ, E. and MASON, D. M. (2004) The law of the iterated logarithm for the integrated squared deviation of a kernel density estimator. *Bernoulli* **10** 721–752 MR2076071

[12] HANSON, D. L. and WRIGHT, F. T. (1971) A bound on the tail probabilities for quadratic forms in independent random variables. *Ann. Math. Statist.* **42** 1079–1083 MR279864

[13] HOEFFDING, W. (1948) A class of statistics with asymptotically normal distribution. *Ann. Math. Statist.* **19** 293–325 MR26294

[14] ITÔ K. (1951) Multiple Wiener integral. *J. Math. Soc. Japan* **3**. 157–164 MR44064

[15] LATAŁA, R. (2005) Estimates of moments and tails of Gaussian chaoses. Available at `http://www.arxiv.org/abs/math.PR/0505313`

[16] LEDOUX, M. (1996) On Talagrand deviation inequalities for product measures. *ESAIM: Probab. Statist.* **1**. 63–87. Available at `http://www.emath./fr/ps/`. MR1399224

[17] MAJOR, P. (1981) Multiple Wiener–Itô integrals. *Lecture Notes in Mathematics* **849**, Springer Verlag, Berlin, Heidelberg, New York, MR611334

[18] MAJOR, P. (2005) An estimate about multiple stochastic integrals with respect to a normalized empirical measure. *Studia Scientarum Mathematicarum Hungarica.* **42** (3) 295–341

[19] MAJOR, P. (2006) An estimate on the maximum of a nice class of stochastic integrals. *Probability Theory and Related Fields.* **134** (3) 489–537

[20] MAJOR, P. (2005) On a multivariate version of Bernstein's inequality. Submitted to *Annals of Probability*. Available at the homepage `http://www.renyi.hu/~major`

[21] MAJOR, P. (2005) A multivariate generalization of Hoeffding's inequality. Submitted to *Annals of Probability*. Available at the homepage `http://www.renyi.hu/~major`

[22] MAJOR, P. (2005) On the tail behaviour of multiple random integrals and degenerate *U*-statistics. (manuscript for a future Lecture Note) Available at the homepage `http://www.renyi.hu/~major`

[23] MAJOR, P. and REJTŐ, L. (1988) Strong embedding of the distribution function under random censorship. *Annals of Statistics*, **16**, 1113–1132 MR959190

[24] MAJOR, P. and REJTŐ, L. (1998) A note on nonparametric estimations. *A volume in Honour of Miklós Csörgő.* North Holland 759–774 MR1661516

[25] MALYSHEV, V. A. and MINLOS, R. A. (991) *Gibbs Random Fields. Method of cluster expansion.* Kluwer, Academic Publishers, Dordrecht MR1191166

[26] MASSART, P. (2000) About the constants in Talagrand's concentration inequalities for empirical processes. *Annals of Probability* **28**, 863–884 MR1782276

[27] MC. KEAN, H. P. (1973) Wiener's theory of non-linear noise. in *Stochastic Differential Equations* SIAM–AMS Proc. 6 197–209 MR353471

[28] POLLARD, D. (1984) *Convergence of Stochastic Processes.* Springer Verlag, New York MR762984

[29] ROTA, G.-C. and WALLSTROM, C. (1997) Stochastic integrals: a combinatorial approach. *Annals of Probability* **25** (3) 1257–1283 MR1457619

[30] STRASSEN, V. (1966) A converse to the law of the iterated logarithm. *Z. Wahrscheinlichkeitstheorie* **4** 265–268 MR200965

[31] SURGAILIS, D. (2003) CLTs for polynomials of linear sequences: Diagram formula with illustrations. in *Long Range Dependence* 111–127 Birkhäuser, Boston, Boston, MA. MR1956046

[32] TAKEMURA, A. (1983) Tensor Analysis of ANOVA decomposition. *J. Amer. Statist. Assoc.* **78**, 894–900 MR727575

[33] TALAGRAND, M. (1994) Sharper bounds for Gaussian and empirical processes. *Ann. Probab.* **22**, 28–76 MR1258865

[34] TALAGRAND, M. (1996) New concentration inequalities in product spaces. *Invent. Math.* **126**, 505–563 MR1419006